\input amstex
\documentstyle{amsppt}
\magnification=\magstep1

\NoBlackBoxes
\TagsAsMath

\pageheight{9.0truein}
\pagewidth{6.5truein}

\long\def\ignore#1\endignore{\par DIAGRAM\par}
\long\def\ignore#1\endignore{#1}

\ignore
\input xy \xyoption{matrix} \xyoption{arrow}
          \xyoption{curve}  \xyoption{frame}
\def\edge{\ar@{-}}
\def\dttdar{\ar@{.>}}

\def\levelpool#1{\save [0,0]+(-3,3);[0,#1]+(3,-3)
                  **\frm<10pt>{.}\restore}
\def\dashedge{\ar@{--}}

\def\dshdar{\ar@{-->}}

\def\dropup#1{\save+<0ex,4ex> \drop{#1} \restore}
\def\dropvert#1#2{\save+<0ex,#1ex>\drop{#2}\restore}

\def\loopNE{\ar@'{@+{[0,0]+(6,2)} @+{[0,0]+(10,10)}
@+{[0,0]+(2,6)}}}
\def\loopNW{\ar@'{@+{[0,0]+(-2,6)} @+{[0,0]+(-10,10)}
@+{[0,0]+(-6,2)}}}
\def\loopSW{\ar@'{@+{[0,0]+(-6,-2)} @+{[0,0]+(-10,-10)}
@+{[0,0]+(-2,-6)}}}
\def\loopSE{\ar@'{@+{[0,0]+(2,-6)} @+{[0,0]+(10,-10)}
@+{[0,0]+(6,-2)}}}

\def\loopNNE{\ar@'{@+{[0,0]+(4,2)} @+{[0,0]+(6,11)}
@+{[0,0]+(0,6)}}}
\def\loopSSW{\ar@'{@+{[0,0]+(-4,-2)} @+{[0,0]+(-6,-12)}
@+{[0,0]+(0,-6)}}}
\def\loopSSE{\ar@'{@+{[0,0]+(0,-6)} @+{[0,0]+(6,-11)}
@+{[0,0]+(4,-2)}}}
\endignore

\def\la{{\Lambda}}
\def\lamod{\Lambda\text{-}\roman{mod}}

\def\AA{{\Bbb A}}
\def\CC{{\Bbb C}}
\def\PP{{\Bbb P}}
\def\SS{{\Bbb S}}
\def\ZZ{{\Bbb Z}}
\def\NN{{\Bbb N}}

\def\hom{\operatorname{Hom}}
\def\aut{\operatorname{Aut}}

\def\top{\operatorname{top}}

\def\ann{\operatorname{ann}}

\def\stab{\operatorname{Stab}}

\def\GL{\operatorname{GL}}

\def\autlap{\operatorname{Aut}_\la(P)}
\def\End{\operatorname{End}}

\def\boldP{{\bold P}}

\def\C{{\Cal C}}

\def\F{{\Cal F}}

\def\T{{\Cal T}}

\def\M{{\Cal M}}
\def\fancyO{{\Cal O}}

\def\S{{\sigma}}

\def\T{{\Cal T}}

\def\bd{\bold {d}}

\def\modlad{\operatorname{\bold{Mod}}_{\bold{d}}(\Lambda)}

\def\toptbd{\operatorname{\bold{Mod}^T_{\bold{d}}}}

\def\grasstd{\operatorname{\frak{Grass}}^T_d}

\def\Gr{\operatorname{Gr}}
\def\grasstbd{\operatorname{\frak{Grass}^T_{\bold{d}}}}
\def\autlap{\aut_\Lambda(P)}

\def\unirad{\bigl(\aut_\la(P)\bigr)_u}
\def\grassS{\operatorname{\frak{Grass}}(\S)}

\def\grassSS{\operatorname{\frak{Grass}}(\SS)}

 \def\biggrass{\operatorname{GRASS}(\la,\bd)}

\def\degen{\le_{\text{deg}}}
\def\underbardim{\operatorname{\underline{dim}}}
\def\id{\operatorname{id}}
\def\term{\operatorname{end}}
\def\start{\operatorname{start}}
\def\Schu{\operatorname{Schu}}

\def\flag{\operatorname{\frak{Flag}}}
\def\Schu{\operatorname{\Schu}}

\def\grasstd{\operatorname{\frak{Grass}}^T_d}
\def\gradgrasstd{\operatorname{grad- \frak{Grass}}^T_d}
\def\gradautlap{\operatorname{grad-Aut}_{\la}(P)}
\def\grasstbd{\operatorname{\frak{Grass}}^T_{\bold d}}
\def\GRASS{\operatorname{GRASS}}

\def\grassd{\GRASS_d(\Lambda)}
\def\grassbd{\GRASS_{\bold d}(\Lambda)}

\def\Gr{\operatorname{Gr}}
\def\boldgrasstd{\operatorname{\frak{Grass}}^T_{\bold{d}}}

\def\maxmoduli{\operatorname{\frak{ModuliMax}^T_{\bold{d}}}}

\def\maxmoduliM{\operatorname{\frak{Maxtopdeg}}}
\def\autlap{\aut_\Lambda(P)}
\def\autlat{\aut_\Lambda(T)}

\def\unirad{\bigl(\autlap \bigr)_u}

\def\grassS{\operatorname{\frak{Grass}}(\S)}

\def\grassSS{\operatorname{\frak{Grass}}(\SS)}

\def\biggrass{\GRASS}

\def\underbardim{\operatorname{\underline{dim}}}
\def\id{\operatorname{id}}
\def\term{\operatorname{end}}
\def\start{\operatorname{start}}
\def\Schu{\operatorname{Schu}}
\def\spec{\operatorname{Spec}}

\def\Adrialebruyn{{\bf 1}}
\def\BHTone{{\bf 2}}
\def\BHTtwo{{\bf 3}}
\def\Bongadv{{\bf 4}}
\def\Bongtrond{{\bf 5}}
\def\BoHZone{{\bf 6}}
\def\BoHZtwo{{\bf 7}}
\def\Bor{{\bf 8}}

\def\CBS{{\bf 9}}
\def\Der{{\bf 10}}
\def\DHZW{{\bf 11}}
\def\DerWey{{\bf 12}}
\def\Geiss{{\bf 13}}

\def\Har{{\bf 14}}
\def\Hille{{\bf 15}}
\def\GeomI{{\bf 16}}
\def\menace{{\bf 17}}
\def\classifying{{\bf 18}}
\def\topstableI{{\bf 19}}
\def\hier{{\bf 20}}

\def\KacI{{\bf 21}}
\def\KacII{{\bf 22}}
\def\King{{\bf 23}}
\def\Kra{{\bf 24}}
\def\MumFogKir{{\bf 25}}
\def\New{{\bf 26}}

\def\PopVin{{\bf 27}}
\def\Rein{{\bf 28}}
\def\Reinsurv{{\bf 29}}
\def\Riedt{{\bf 30}}
\def\Rosone{{\bf 31}}
\def\Rostwo{{\bf 32}}
\def\Scho{{\bf 33}}
\def\Schro{{\bf 34}}
\def\Zwara{{\bf 35}}

\topmatter

\title Fine and coarse moduli spaces in the representation theory of finite dimensional algebras \endtitle

\rightheadtext{moduli spaces}

\author B. Huisgen-Zimmermann
\endauthor

\dedicatory Dedicated to Ragnar-Olaf Buchweitz on the occasion of his seventieth birthday  \enddedicatory

\thanks The author was partially supported by a grant from the National Science Foundation \endthanks

\address Department of Mathematics, University of California, Santa
Barbara, CA 93106-3080 \endaddress

\abstract We discuss the concepts of fine and coarse moduli spaces in the context of finite dimensional algebras over algebraically closed fields.  In particular, our formulation of a moduli problem and its potential strong or weak solution is adapted to classification problems arising in the representation theory of such algebras. We then outline and illustrate a dichotomy of strategies for concrete applications of these ideas.  One method is based on the classical affine variety of representations of fixed dimension, the other on a projective variety parametrizing the same isomorphism classes of modules. We state sample results and give numerous examples to exhibit pros and cons of the two lines of approach.  The juxtaposition highlights differences in techniques and attainable goals. 
\endabstract
   
\endtopmatter

\document

\head 1. Introduction and notation \endhead

The desire to describe/classify the objects of various algebro-geometric categories via collections of invariants is a red thread that can be traced throughout mathematics.  Prominent examples are the classification of  similarity classes of matrices in terms of normal forms, the classification of finitely generated abelian groups in terms of annihilators of their indecomposable direct summands, and the classification of varieties of fixed genus and dimension up to isomorphism or birational equivalence,  etc., etc. --  the reader will readily extend the list.   In each setting, one selects an equivalence relation on the collection of objects to be sorted;  the ``invariants" one uses to describe the objects are quantities not depending on the choice of representatives from the considered equivalence classes; and the chosen data combine to finite parcels that identify these classes, preferably without redundancy.  In case the relevant parcels of invariants consist of discrete data  --  as in the classification of finitely generated abelian groups up to isomorphism for instance  --   there is typically no need for additional tools to organize them.  By contrast, if the objects to be classified involve a base field $K$ and their invariants are structure constants residing in this field  --  suppose one has established a one-to-one correspondence between the equivalence classes of objects and certain points in an affine or projective space over  $K$  --  it is natural to ask whether these invariants trace an algebraic variety over $K$.  In the positive case, one is led to an analysis of how structural properties of the objects under discussion are mirrored by geometric properties of the pertinent parametrizing variety.  The least one hopes for is some form of ``continuous" bijective dependence of the equivalence classes of objects on their classifying parameters (what one means by ``continuity" is made precise through the notion of a ``family"), preferably satisfying a universal property.    

The prototypical example of a highly successful classification of this ilk goes back to Riemann:   In 1857, he classified the isomorphism classes of
nonsingular complex projective curves of fixed genus $g
\ge 2$ in terms of what he called ``moduli".  According to Riemann, the
isomorphism class of a nonsingular curve of genus $g$ ``haengt von $3g
- 3$ stetig veraenderlichen Groessen ab, welche die Moduln dieser
Klasse genannt werden sollen".  It took about a century for the term ``moduli" to be given a precise meaning, namely, as elements
of a fine or coarse moduli space.  Such moduli spaces were
axiomatically introduced by Mumford in the 1960s.   At the beginning of Section 2, we will roughly sketch the underlying idea, adapted to the representation theory of a finite dimensional algebra $\la$, in order to motivate a first round of tool-building.  Precise definitions of moduli spaces are given in Section 4, which follows the transparent exposition of Newstead \cite{\New}. 

To delineate our goals: Our present interest is in an improved understanding of select classes of representations of a basic finite dimensional algebra $\la$ over an algebraically closed field $K$ by way of moduli.  We may assume without loss of generality that $\la = KQ/I$, where $Q$ is a quiver with vertex set $Q_0 = \{e_1, \dots, e_n\}$, and $I$ is an admissible ideal in the path algebra $KQ$.   Our primary objective here is to review and compare presently available techniques and results that harness fine or coarse moduli spaces for the classification of finite dimensional representations of such an algebra $\la$.  A discussion from a unified perspective should, in particular, make the subject more accessible to newcomers to the area;  to meet this purpose, we will include some elementary observations, to be skipped by experts.  A secondary aim is to promote a series of problems which appear to be ``next in line" towards broadening the impact of the general moduli machinery on representation theory.  Some of these problems  --  those aiming at ``generic classification"  --  extend an investigation that was initiated by Kac in the early 1980s (see \cite{\KacI, \KacII}) and picked up by Schofield \cite{\Scho}, Derksen-Weyman \cite{\DerWey}, Reineke \cite{\Reinsurv} and others in the case $I = 0$, by Schr\"oer \cite{\Schro}, Crawley-Boevey-Schr\"oer \cite{\CBS},  Babson-Thomas and the author \cite{\BHTtwo} in the general situation.  

Here is an outline of the article:  In Section 1, we revisit two starting points for a geometric classification of finite dimensional $\la$-modules.  We first review the classical affine variety $\modlad$ parametrizing the (left) $\la$-modules with fixed dimension vector $\bd$ (we refer to it as Parametrization A); next we turn to the projective variety  $\grassbd$ parametrizing the same isomorphism classes of modules (Parametrization B).   In either case, the parametrizing variety comes equipped with an algebraic group action, the orbits of which are in bijective correspondence with the isomorphism classes of modules under consideration.  However, the widely different structures of these varieties and their respective acting groups give the two points of departure distinct types of potential, on some occasions yielding alternate roads to the same conclusion.  In both settings, one observes that the group action can hardly ever be factored out of the original parametrizing variety in a geometrically meaningful manner, which prompts us to include a brief general discussion of quotients of algebraic varieties modulo actions of algebraic groups in Section 3.   This section overlaps with expository articles by Bongartz \cite{\Bongtrond} and Geiss \cite{\Geiss}.  The modest overlap is required for a consistent development of the subsequent ideas.

Then we return to Riemann's classification program and discuss/exemplify the concepts of a fine/coarse moduli space in the representation-theoretic context (Section 4).  To date, there are two different strategies to get mileage out of the conceptual framework.  In light of the fact that fine or coarse moduli spaces for the full collection of isomorphism classes of $\la$-modules with a given dimension vector hardly ever exist, each method proposes a mode of slicing $\lamod$ so as to extract portions on which the conceptual vehicle of moduli spaces acquires traction.  The strategies of slicing take advantage of the particulars of the initial parametrizing setups, and hence, in each case, specific methodology is called for to match the target.  Since there exist two prior survey articles dealing with Approach A, by Geiss \cite{\Geiss} and Reineke \cite{\Reinsurv}, we will give more room to Approach B in the present overview. 

One of the methods mimicks a strategy Mumford used in the classification of vector bundles on certain projective varieties.  It was adapted  to the representation-theoretic setting by King in \cite{\King} (see Section 5).  Starting with an additive function $\theta: K_0(\la) = \ZZ^n \rightarrow \ZZ$, King focuses on the $\la$-modules with dimension vector $\bd$ which are $\theta$-{\it stable\/}, resp\. $\theta$-{\it semistable\/}; interprets these stability conditions in terms of the behavior of $\theta$ on submodule lattices; and shows how to apply techniques from geometric invariant theory to secure a fine, resp. coarse, moduli space for $\theta$-(semi)stable modules.  The resulting stability classes are not a priori representation-theoretically distinguished, whence a fundamental challenge lies in ``good" choices of the function $\theta$ and a solid grasp of the corresponding $\theta$-(semi)stable modules.   As this method is based on the affine parametrizing variety $\modlad$, crucially leaning on the features of this setup, it will be labeled Approach A.  So far, its main applications are to the hereditary case $\la = KQ$, even though, in principle, King extended the method to include arbitrary path algebras modulo relations.

By contrast, the second approach (labeled Approach B and described in Sections 6-8) starts with classes $\C$ of modules over $\la = KQ/I$ which are cut out by purely representation-theoretic features, and aims at understanding these classes through an analysis of the subvarieties of $\grassbd$ that encode them.  The name of the game is to   exploit projectivity of the parametrizing variety and the typically large unipotent radical of the acting group to find useful necessary and sufficient conditions for the existence of a geometric quotient of the subvariety encoding $\C$, and to subsequently establish such a quotient as a moduli space that classifies the representations in $\C$ up to isomorphism.  Simultaneously, one seeks theoretical and/or algorithmic access to moduli spaces whenever existence is guaranteed.

In describing either method, we state sample theorems witnessing viability and illustrate them with examples.  Each of the two outlines will conclude with a discussion of pros and cons of the exhibited approach.  

{\it Acknowledgements.}  I wish to thank the organizers of the Auslander Conference and Distinguished Lecture Series (Woods Hole, April 2012), K.~ Igusa, A.~ Martsinkovsky, and G.~ Todorov, and the organizers F.~ Bleher and C.~ Chindris of the Conference on Geometric Methods in Representation Theory (University of Missouri-Columbia, November 2012) for having provided me with congenial venues for the expository lectures that gave rise to these notes.  
\bigskip

{\it Further conventions.}  Throughout, $\la$ will be a basic finite dimensional algebra over an algebraically closed field $K$, and $J$ will denote the Jacobson radical of $\la$. We thus do not lose any generality in assuming that $\la = KQ/I$ for a quiver $Q$ and an admissible ideal $I$ of the path algebra $KQ$.  The vertex set $Q_0 = \{ e_1, \dots, e_n\}$ of $Q$ will be identified with a full set of primitive idempotents of  $\la$.  Moreover, we let $S_i = \la e_i / Je_i$ be the corresponding representatives of the simple modules. The {\it absolute value\/} of a dimension vector $\bd = (d_1, \dots, d_n)$  is $|\bd| = \sum_i d_i$.   

We will systematically identify isomorphic semisimple modules.  The {\it top\/} of a (left) $\la$-module $M$ is $\top(M) = M / JM$.  The {\it radical layering\/} of $M$ is the sequence of semisimple modules $\SS(M) = \bigl(J^l M/ J^{l+1}M \bigr)_{0 \le l \le L}$, where $L+1$ is the Loewy length of $\la$.  In particular, the zero-th entry of $\SS(M)$ equals the top of $M$.

For our present purpose, it suffices to consider classical quasi-projective varieties.
By a {\it subvariety\/} of a such a variety we will mean a locally closed subset.

\head 2.  Affine and projective parametrizations of the $\la$-modules of dimension vector $\bd$ \endhead

Suppose that $\C$ is a class of objects in some algebro-geometric category, and let $\sim$ be an equivalence relation on $\C$. 

\definition{Riemann's classification philosophy in loose terms}  
\smallskip

\noindent {\bf (I)}  Identify discrete invariants of the objects in $\C$, in order to subdivide $\C$ into finitely many (or countably many) subclasses $\C_i$, the objects of which are sufficiently akin to each other to allow for a normal form characterizing them up to the chosen equivalence.  
\smallskip

\noindent {\bf (II)}  For each index $i$, find an algebraic variety $V_i$, together with a bijection 
$$V_i \ \ \longleftrightarrow \ \ \{\text{equivalence classes in\ \,} \C_i \},$$
which yields a {\it continuous parametrization\/} of the equivalence classes of objects in $\C_i$.  (The idea of ``continuity" will be clarified in Section 4.  Typically, such a parametrization will  --  a priori or a posteriori  --  be a classification of normal forms.)  Once a parametrization that meets these ciriteria is available, explore potential {\it universal properties\/}.  Moreover, investigate the interplay between the geometry of $V_i$ on one hand and structural properties of the modules in $\C_i$ on the other.
\enddefinition 
  
We will focus on the situation where $\C$ is a class of representations of $\la$.  In this situation, the most obvious equivalence relation is isomorphism, or graded isomorphism if applicable. Riemann's philosophy then suggests the following as a first step:  Namely, to tentatively parametrize the isomorphism classes of modules with fixed dimension vector in {\it some\/} plausible way by a variety.  We will review two such parametrizations, both highly redundant in the sense that large subvarieties map to single isomorphism classes in general.  In each case, the considered parametrizing variety carries a morphic action by an algebraic group $G$ whose orbits capture the redundancy; in other words, the $G$-orbits are precisely the sets of points indexing objects from the same isomorphism class of modules.   Since each of these settings will have advantages and downsides compared with the other, it will be desirable to shift data back and forth between them.  Such a transfer of information between Scenarios A and B will turn out to be optimally smooth.  We will defer a detailed discussion of this point to the end of Section 3, however, since we wish to specifically address the passage of information concerning quotients by the respective group actions.     

\smallskip

\noindent {\bf (A) The classical affine parametrization of the isomorphism classes of $\la$-modules with dimension vector $\bd$} 
\smallskip

This setup is well-known and much-used.  To our knowledge, the first prominent application to the representation theory of finite dimensional algebras was in the proof of Gabriel's Theorem pinning down the path algebras of finite representation type.  

 \definition{The affine parameter variety and its group action}
 \smallskip

(1)  Let $Q_1$ be the set of arrows of $Q$, and let
\smallskip

$\modlad =$
\smallskip

\qquad \ \ $ \{ x = (x_{\alpha})_{\alpha \in Q_1} \mid  \text{the\ } x_\alpha \in M_{d_{\term(\alpha)}\times d_{\start(\alpha)}}(K)  \text{\ satisfy the relations in\ } I\}$.
\smallskip

\noindent  Here $M_{r \times s}(K)$ denotes the space of $r \times s$ matrices over $K$.
\smallskip

(2) The group action:  Set $\GL(\bd) = \prod_{1 \le i \le n} \GL_{d_i}(K)$, and consider the following action of $\GL(\bd)$ on $\modlad$:  For $g = (g_1, \dots, g_n) \in \GL(\bd)$ and $x  = (x_{\alpha})\in \modlad$, define
\smallskip 
\centerline{$g.x = \bigl(g_{\term(\alpha)} x_\alpha g_{\start(\alpha)}^{-1}\bigr)_{\alpha \in Q_1}$.}
\enddefinition

Evidently, $\modlad$ is a Zariski-closed subset of the affine $K$-space of dimension $\sum_{\alpha \in Q_1} d_{\start(\alpha)} d_{\term(\alpha)}$, the points of which determine $\la$-module structures on the vector space $K^{|\bd|} = \bigoplus_{1 \le i \le n} K^{d_i}$ via $\alpha \bigl(\sum_{1 \le i \le n} v_i \bigr) = x_\alpha \bigl(v_{\start(\alpha)} \bigr) \in K^{d_{\term(\alpha)}}$ for any arrow $\alpha$ and $v_i \in K^{d_i}$.  Clearly, the fibers of the resulting map from $\modlad$ to the set of isomorphism classes of modules with dimension vector $\bd$ are precisely the orbits of the described $\GL(\bd)$-action on $\modlad$.  Thus, we obtain a one-to-one correspondence between the $\GL(\bd)$-orbits of $\modlad$ on one hand and the isomorphism classes of $\la$-modules with dimension vector $\bd$ on the other.  Moreover, we observe that the considered group action is morphic, meaning that the pertinent map $\GL(\bd) \times \modlad \rightarrow \modlad$ is a morphism of varieties.
\medskip

\noindent {\bf (B) The projective parametrization of the same set of isomorphism classes} 
\smallskip 

An alternate parametrizing variety for the same isomorphism classes of modules was introduced by Bongartz and the author in \cite{\BoHZone, \BoHZtwo}, together with a morphic algebraic group action whose orbits, in turn, are in one-to-one correspondence with these isomorphism classes.

\definition{The projective parameter variety and its group action} 
\smallskip 

(1)  Let $\bold P = \bigoplus_{1 \le i \le n} (\la e_i)^{d_i}$ (the smallest projective $\la$-module admitting arbitrary modules with dimension vector $\bd$ as quotients modulo suitable submodules), and define
$$\grassbd = \{ C \in \Gr(d', \, \boldP) \mid \  _{\la}C \subseteq  {}_{\la}\boldP \text{\ with\ } \underline{\dim}\, \boldP /C = \bd \},$$
where $d' =  \dim \boldP-|\bd|$ and $\Gr(d', \, \boldP)$ is the Grassmann variety of all $d'$-dimensional subspaces of the $K$-vector space $\boldP$.
\smallskip

(2)  The group action:  Let $\aut_\la(\boldP)$ be the automorphism group of $\boldP$, and consider the canonical action on $\grassbd$ given by $f.C = f(C)$.
\enddefinition

This time, we are looking at a Zariski-closed subset of the classical Grassmann variety $\Gr(d', \, \boldP)$; in particular, $\grassbd$ is a projective variety.  Again, we have an obvious map from the variety $\grassbd$ to the set of isomorphism classes of $\la$-modules with dimension vector $\bd$, namely $\rho: C \mapsto [\boldP/C]$.  By the choice of $\boldP$, every module $M$ with dimension vector $\bd$ is of the form $M \cong \boldP / C$ for some point $C \in \grassbd$. Moreover, the fibers of $\rho$ again coincide with the orbits of the group action; indeed, two modules $\boldP/C$ and $\boldP/D$ are isomorphic iff $C$ and $D$ belong to the same orbit, this time the $\aut_\la(\boldP)$-orbit of $\grassbd$.  Moreover, the group action is in turn morphic.

Recall that the {\it unipotent radical\/} of a linear algebraic group is the unique largest normal connected unipotent subgroup.  The group is called {\it reductive\/} if its unipotent radical is trivial. 
In contrast to the reductive group $\GL(\bd)$ acting in the affine case, the linear group $\aut_\la(\boldP)$ has a large unipotent radical in most interesting cases.  Namely, the unipotent radical, $\aut_\la(\boldP)_u$, equals the subgroup $\{\id + h \mid  h \in \hom_\la(\boldP, J \boldP)\}$.  We observe moreover that $\aut_\la(\boldP) \cong \GL(\boldP  / J\boldP) \ltimes  \aut_\la(\boldP)_u$. 
\bigskip

\head 3.  Quotient varieties on the geometric market  --  generalities and representation-theoretic particulars\endhead

In Section 2, we have, in both cases, arrived at a scenario that is frequently encountered in connection with classification problems:  One starts with a collection of algebro-geometric objects which one wishes to classify up to an equivalence relation --  in our case the objects are representations with fixed dimension vector and the preferred equivalence relation is isomorphism.  On the road, one arrives at a setup that places the equivalence classes of objects into a natural one-to-one correspondence with the orbits of an algebraic group action on a parametrizing variety.   Such a  scenario, of course, triggers the impulse to factor the group action out of the considered variety.   To say it in different words:  The idea is to reduce the orbits of the group action to points in a new variety which is related to the original one by a universal property which takes the geometry into account. 

The crux lies in the fact that the topological quotient of $\modlad$ modulo $\GL(\bd)$, (resp.~ of  $\grassbd$ modulo $\aut_\la(\boldP)$), relative to the Zariski topology, hardly ever carries a variety structure, at least not one that merits the label ``quotient variety".  To cope with this difficulty in a broad spectrum of situations, algebraic geometers introduced quotients of various levels of stringency.  Not surprisingly, the underlying guideline is this:  The closer the Zariski topology of a ``quotient variety" comes to that of the topological quotient, the better. We will touch this subject only briefly and refer the reader to the  survey by Bongartz \cite{\Bongtrond} and the exposition by Popov and Vinberg \cite{\PopVin}.

\definition{Categorical and geometric quotients}  Let $X$ be an algebraic variety, and let $G$ be a linear algebraic group acting morphically on $X$.  
\smallskip

\noindent {\bf (1)}  A {\it categorical {\rm{(}}or algebraic{\rm{)}}} quotient of $X$ by $G$ is a morphism $\psi: X \rightarrow  Z$ of varieties such that $\psi$ is constant on the orbits of $G$, and every morphism $\psi':  X \rightarrow Y$ which is constant on the $G$-orbits factors uniquely through $\psi$.  Write $Z = X/\!/G$ in case such a quotient exists.
\smallskip

\noindent {\bf (2)}  A categorical quotient $\psi: X \rightarrow X/\!/G$ is called an {\it orbit space\/} for the action in case the fibers of $\psi$ coincide with the orbits of $G$ in $X$.
\smallskip 

\noindent {\bf (3)}  A  {\it geometric quotient of $X$ by $G$\/} is an open surjective morphism $\psi:  X \rightarrow Z$, whose fibers equal the orbits of $G$ in $X$, such that, moreover, for every open subset $U$ of $Z$, the comorphism $\psi^\circ$ induces an algebra isomorphism from the ring $\fancyO(U)$ of regular functions on $U$ to the ring $\fancyO \bigl( \psi^{-1}(U) \bigr)^G$ of $G$-invariant regular functions on $\psi^{-1}(U)$.
\enddefinition

It is easy to see that a geometric quotient is an orbit space, and hence, in particular, is a categorical quotient.  This guarantees uniqueness in case of existence.  We give two elementary examples in order to build intuition:  For $n \ge 2$, the conjugation action of $\GL_n(K)$ on the variety of $n \times n$ matrices has a categorical quotient, which, however, fails to be an orbit space.   Given a linear algebraic group $G$ and any closed subgroup $H$, the right translation action of $H$ on $G$ has a geometric quotient;  in particular, the points of this quotient may be identified with the left cosets of $H$ in $G$.

One readily verifies that the Zariski topology on a geometric quotient coincides with the quotient topology.  So, in light of the above guideline, existence of a geometric quotient is the best possible outcome whenever we look for a quotient of a subvariety of $\modlad$ modulo $GL(\bd)$ or of a subvariety of $\grassbd$ modulo $\aut_\la(\boldP)$.  On the other hand, an orbit space for a suitable action-stable subvariety is the least we require in order to implement Riemann's idea.  Evidently, 
\smallskip

\centerline{$\bullet$ \it the existence of an orbit space implies closedness of all orbits,}
\smallskip
\noindent which places a strong necessary condition on potential scenarios of success.

Let us take a look at our two parametrizations of the $\la$-modules with dimension vector $\bd$.  Here is what Geometric Invariant Theory grants us in the affine setting:
Namely, every morphic action of a reductive linear algebraic group $G$ on an affine variety  $X$ has a categorical quotient (see, e.g., \cite{\New, Chapter 3}).  The pivotal asset of this setup lies in the fact that the ring $K[X]^G$ of $G$-invariant regular functions (i.e., of regular functions $f: X \rightarrow K$ such that f(gx) = f(x), for all $g \in G$ and $x \in X$) is finitely generated over $K$.  We will repeatedly refer to this result. 

\proclaim{Theorem 3.1} {\rm{[Haboush, Hilbert, Mumford, Nagata, Weyl, et al.]}}  Suppose that $X$ is an affine variety with coordinate ring $K[X]$and $G$ a reductive group acting morphically on $X$.  Then the canonical map 
$$\psi: \spec K[X]  \rightarrow \spec K[X]^G$$ 
is a categorical quotient $X/\!/G$.  Moreover, the points of $X/\!/G$ are in one-to-one correspondence with the closed $G$-orbits of $X$.
\endproclaim

 In particular, Theorem 3.1 guarantees a categorical quotient $\modlad /\!/ \GL(\bd)$.  At first glance, this conclusion may look better than it is, since the only closed orbit in $\modlad$ is that of the semisimple module of dimension vector $\bd$. Indeed, given any module $M$ and any submodule $N \subseteq M$, the $\GL(\bd)$-orbit corresponding to the direct sum $N \oplus M/N$ in $\modlad$ is contained in the closure of the $\GL(\bd)$-orbit corresponding to $M$.  So, by the theorem, $\modlad /\!/ \GL(\bd)$ is a singleton.
Expressed differently:  The catch lies in the fact that the ring of $\GL(\bd)$-invariant regular functions on $\modlad$ equals the field $K$ of constants, and hence has only a single prime ideal.  The response of Geometric Invariant Theory to such a sparsity of closed orbits is to pare down the parametrizing variety and, in tandem, to relax the invariance requirements placed on the regular functions that are expected to separate the orbits, so as to obtain a larger algebra of functions that may be used to construct a useful quotient. 
\smallskip  

In order to benefit from the fact that different arsenals of techniques apply to our two parametrizations, we first explain how to move back and forth between them.

  \proclaim{Proposition 3.2. Information transfer between Parametrizations A and B} {\rm (see \cite{\BoHZtwo, Proposition C})}
\smallskip

Consider the one-to-one correspondence between the orbits of $\grassbd$ on one hand and $\modlad$ on the other, which
assigns to any orbit $\aut_\la(\boldP).C \subseteq \grassbd$ the orbit
$\GL(\bd).x \subseteq \modlad $  representing the same
$\la$-module up to isomorphism.  This correspondence extends to an inclusion-pre\-serv\-ing bijection
$$\Phi: \{ \aut_\la(\boldP) \text{-stable subsets of\ } \grassbd \} \rightarrow
\{\GL(\bd)\text{-stable subsets of\ } \modlad \}$$  
which preserves and reflects
openness, closures, connectedness, irreducibility, and types of
singularities.  

Moreover, let $X$ be a $\GL(\bd)$-stable subvariety of $\modlad$, with corresponding \linebreak $\aut_\la(\boldP)$-stable subvariety $\Phi(X)$ of $\grassbd$.  Then $X$ has an algebraic  quotient {\rm{(}}resp., orbit space/geometric quotient{\rm{)}} by $\GL(\bd)$ if and only if $\Phi(X)$ has an algebraic quotient {\rm{(}}resp., orbit space/geometric quotient{\rm{)}} by $\aut_\la(\boldP)$.   In case of existence, the quotients are isomorphic and have the same separation properties relative to action stable subvarieties of $X$ and $\Phi(X)$, respectively. 
\endproclaim

The transfer result thus allows us to symmetrize the unhelpful conclusion we drew from Theorem 3.1.  The projective variety $\grassbd$ has a categorical quotient by $\aut_\la(\boldP)$, and this quotient is isomorphic to $\modlad /\!/ \GL(\bd)$, a singleton.

Where should we go from here?  We are on the outlook for interesting subvarieties of $\modlad$, resp\. $\grassbd$ which are stable under the pertinent group actions and have the property that all orbits are relatively closed.  Proposition 3.2  tells us that we may interchangeably use the two settings,  A and B, in this quest.

In Sections 5 and 6, 7 we will review and illustrate two different methods to identify subvarieties of this ilk.  But first we will flesh out the vague classification philosophy presented in Section 1.

\head 4. Rendering Riemann's classification philosophy more concrete \endhead

The current understanding of Riemann's ``moduli" views them as ``elements of a fine or coarse moduli space".  The two notions of moduli space, one significantly stronger than the other, were introduced and put to use by Mumford in the 1960's  (see the standard GIT text \cite{\MumFogKir}).   We will follow Newstead's exposition \cite{\New}. 

Both types of moduli spaces build on the concept of a {\it family of objects parametrized by an algebraic variety\/}.  The upcoming definition clarifies the idea of a {\it continuous parametrization\/}, as opposed to a random indexing of objects by the points of a variety.  The (only) plausible definition of a {\it family\/} in the representation-theoretic context was put forth by King in \cite{\King}.  

\definition{Definition: Families of representations} 

\noindent Set $d = |\bd|$, and let $\sim$ be an equivalence relation on the class of $d$-dimensional $\la$-modules.
\smallskip

\noindent {\bf (1)} A {\it family of $d$-dimensional $\la$-modules parametrized by a variety $X$\/} is a pair $(\Delta, \delta)$, where $\Delta$ is a vector bundle of rank $d$ over $X$, and $\delta$ a $K$-algebra homomorphism $\la \rightarrow \End(\Delta)$.
\smallskip

\noindent {\bf (2)}  Extending $\sim$ to families:  Two such families $(\Delta_1, \delta_1)$ and $(\Delta_2, \delta_2)$ parametrized by the same variety $X$ will be called {\it similar\/} in case, for each $x \in X$, the fibers of $\Delta_1$ and $\Delta_2$ over $x$ are $\sim$ equivalent as $\la$-modules under the structures induced by $\delta_1$ and $\delta_2$, respectively.  We write  $(\Delta_1, \delta_1) \sim \Delta_2, \delta_2)$ in this situation.
\smallskip

\noindent {\bf (3)} Induced families:  Given a family $(\Delta, \delta)$ parametrized by $X$ as above, together with a morphism $\tau: Y \rightarrow X$ of varieties, the pull-back bundle of $\Delta$ along $\tau$ is a family of $\la$-modules parametrized by $Y$ (see the remark below).  It is called the {\it family induced from $(\Delta, \delta)$ by $\tau$\/} and is denoted by $\tau^*(\Delta, \delta)$.
\enddefinition

Here, the  vector bundles considered are what Hartshorne \cite{\Har} calls {\it geometric vector bundles\/}:  This means that $\Delta$ carries the structure of  a variety, and all of the occurring maps -- the bundle projection, the local sections responsible for local triviality, and the compatibility maps for the trivialized patches -- are morphisms of varieties.  The requirement that $\delta(\lambda)$,  for $\lambda \in \la$, be an endomorphism of $\Delta$ just means that $\delta(\lambda): \Delta \rightarrow \Delta$ is a morphism of varieties that respects the fibers of the bundle under the projection map; so we find that  each fiber is indeed endowed with a $\la$-module structure.  Since each $\delta(\lambda)$ is a global morphism from $\Delta$ to $\Delta$, this means that the $\la$-module structures on the individual fibers are compatible in a strong geometric sense, thus justifying the interpretation as a continuous array of modules.  
\smallskip

{\it Remark concerning the pull-back construction}:  Using the corresponding trivializations, we readily check that, for $y \in Y$, the pullback diagram 
$$\xymatrixrowsep{2pc}\xymatrixcolsep{4pc}
\xymatrix{
\Delta \ar[r]^{\pi} &X  \\
\tau^*(\Delta) \ar[r]^{\pi^*} \ar[u] &Y \ar[u]_{\tau}
}$$
\noindent permits us to pull back the $\la$-module structure (stemming from $\delta$) on the fibre $\pi^{-1}(\tau(y))$ of $\Delta$ to a $\la$-module structure on the fiber $(\pi^*)^{-1}(y)$ of $\tau^*(\Delta)$; one verifies that these module structures on the individual fibers of $\tau^*(\Delta)$ are compatible, so as to yield a $K$-algebra homomorphism $\delta^*: \la \rightarrow \End \bigl( \tau^*(\Delta) \bigr)$ that induces them. Set $\tau^*(\Delta, \delta) = (\tau^*(\Delta),\delta^*)$.
\smallskip 

It is easily verified that the definitions of ``family" and ``induced family" satisfy  the functorial conditions spelled out as prerequisites for a well-defined ``moduli problem" in ([\New, Conditions 1.4, p.19]).   Namely:  $\bullet$  The equivalence relation on families boils down to the initial equivalence relation $\sim$ on the target class $\C$, if one identifies a family parametrized by a single point with the corresponding module; in fact, the equivalence relation we introduced under (2) above is the coarsest with this property.  (It is not the most natural option, but the easiest to work with in our context.)   $\bullet$  If $\tau: Y \rightarrow X$ and $\sigma: Z \rightarrow Y$ are morphisms of varieties and $(\Delta, \delta)$ is a family of modules over $X$, then $(\tau \circ \sigma)^*(\Delta, \delta) = \sigma^*\bigl(\tau^*(\Delta, \delta)\bigr)$; moreover $(\id_X)^*$ is the identity on families parametrized by $X$.  $\bullet$ Similarity of families is compatible with the pullback operation, that is: If $(\Delta_1, \delta_1)$ and $(\Delta_2, \delta_2)$ are families parametrized by $X$ with $(\Delta_1, \delta_1) \sim(\Delta_2, \delta_2)$ and $\tau$ is as above, then $\tau^*(\Delta_1,\delta_1) \sim \tau^*(\Delta_2, \delta_2)$. 

\definition{Example 4.1}  Let $\la$ be the Kronecker algebra, i.e., $\la = KQ$, where $Q$ is the quiver 
$\xymatrixcolsep{4pc}\xymatrix{
1 \ar@/^/[r]^{\alpha_1} \ar@/_/[r]_{\alpha_2} &2
}$
,  and take $\bd = (1,1)$.  The non-semisimple $2$-dimensional $\la$-modules form a family indexed by the projective line over $K$.  It can informally be presented as $\bigl(M_{[c_1:c_2]}\bigr)_{[c_1:c_2] \in \PP^1}$ with $M_{[c_1:c_2]} = \la e_1 / \la (c_1 \alpha_1 - c_2 \alpha_2)$.   For a formal rendering in the sense of the above definition, consider the two standard affine patches, $U_j = \{ [c_1:c_2] \in \PP^1 \mid c_j \ne 0\}$, and let $\Delta_j = U_j \times K^2$ for $j = 1,2$ be the corresponding trivial bundles.  To make $\Delta_1$ into a family of $\la$-modules, let $\delta_1: \la \rightarrow \End(\Delta_1)$ be such that $\delta_1 (\alpha_1)$ acts on the fibre above $[c_1:c_2]$ via the matrix $\pmatrix 0&0\\ c_2/c_1&0 \endpmatrix$ and $\delta_1(\alpha_2)$ acts via $\pmatrix 0&0\\ 1&0 \endpmatrix$.  Define $\delta_2: \la \rightarrow \End(\Delta_2)$ symmetrically, and glue the two trivial bundles to a bundle $\Delta$ over $\PP^1$ via the morphism 
$$U_1 \cap U_2 \rightarrow \GL_2(K), \ \ [c_1: c_2] \mapsto \pmatrix c_1/c_2&0\\ 0&c_2/c_1 \endpmatrix.$$  
Observe that the $\delta_j$ are compatible with the gluing, that is, they yield a $K$-algebra homomorphism  $\delta: \la \rightarrow \End(\Delta)$, and thus a family $(\Delta, \delta)$.
\enddefinition

\definition{Definition of fine and coarse moduli spaces} We fix a dimension vector $\bd$,  set $d = |\bd|$, and let $\C$ be a class of $\la$-modules with dimension vector $\bd$.  Denoting by $\C(\modlad)$, resp. $\C(\grassbd)$, the union of all orbits in $\modlad$, resp\. in $\grassbd$, which correspond to the isomorphism classes in $\C$, we assume that $\C(\modlad)$ is a subvariety of $\modlad$ (or, equivalently, that  $\C(\grassbd)$ is a subvariety of $\grassbd$).  Again, we let $\sim$ be an equivalence relation on $\C$ and extend the relation $\sim$ to families as spelled out in the preceding definition.
\smallskip

\noindent {\rm \bf{(1)}}  A {\it fine moduli space\/} classifying $\C$ up to $\sim$ is a variety $X$ with the property that there exists a family $(\Delta, \delta)$ of modules from $\C$ which is parametrized by $X$ and has the following universal property:  Whenever $(\Gamma, \gamma)$ is a family of modules from $\C$ indexed by a variety $Y$, there exists a unique morphism $\tau:Y \rightarrow X$ such that $(\Gamma, \gamma) \sim \tau^* (\Delta, \delta)$.

In this situation, we call $(\Delta,\delta)$ a {\it universal family\/} for our classification problem. 
(Clearly, such a universal family is unique up to $\sim$ whenever it exists.)
\smallskip

\noindent {\rm \bf{(2)}} Specializing to the case where $\sim$ is ``isomorphism" (for the moment), we say that a variety $X$ is a {\it coarse moduli space\/} for the classification of $\C$ up to isomorphism in case $X$ is an orbit space for $\C(\modlad)$ under the $\GL(\bd)$-action (or, equivalently, for $\C(\grassbd)$ under the $\aut_\la(\boldP)$-action).  
\enddefinition

In Section 6, we will also look for moduli spaces classifying classes $\C$ of graded modules up to graded isomorphism.  By this we will mean an orbit space of $\C(\grassd)$ relative to the action of the group of graded automorphisms in $\aut_\la(\bold P)$.  

\definition{Comments 4.2}  Rather than giving the original functorial definitions of fine/coarse moduli spaces, we have introduced these concepts via equivalent characterizations of higher intuitive appeal.
\smallskip

\noindent \noindent {\rm \bf{(i)}} The standard functorial definitions of a fine/coarse moduli space are as follows (cf\. \cite{\New}):  
 
 Consider the contravariant functor 
$$F:  \text{Var = category of varieties over} \ K\ \ \  \longrightarrow\ \ \ \text{category of sets},$$
$$Y \mapsto \{\text{equivalence classes of families of objects from\ } \C \text{\ parametrized by\ } Y\}.$$ 

This functor is representable in the form $F \cong \hom_{\text{Var}}(-, X)$ precisely when    
$X$ is a fine moduli space for our problem.

That a variety $X$ be a coarse moduli space for our problem amounts to the following condition:  There exists a natural transformation $\Phi: F \rightarrow \hom_{\text{Var}}(-, X)$ such that any natural transformation $F \rightarrow \hom_{\text{Var}}(-, Y)$ for some variety $Y$ factors uniquely through $\Phi$. 
\smallskip

\noindent \noindent {\rm \bf{(ii)}}  Our definition of a coarse moduli space $X$ is equivalent to Mumford's in the situations on which we are focusing, but not in general.  We are chipping in the fact that the modules from $\C$ belong to a family $(\Delta, \delta)$ that enjoys the {\it local universal property\/} in the sense of \cite{\New, Proposition 2.13}; indeed, we only need to restrict the tautological bundle on $\modlad$ to $\C(\modlad)$ .  If $X$ denotes the parametrizing variety of $\Delta$, this condition postulates the following: For any family $(\Gamma, \gamma)$ of modules from $\C$, parametrized by a variety $Y$ say, and any $y \in Y$, there is a neighborhood $N(y)$ of $y$ such that the restricted family $\Gamma |_{N(y)}$ is induced from $\Delta$ by a morphism $N(y) \rightarrow X$.  Note that local universality carries no uniqueness requirement.

In classifying graded representations of a graded algebra $\la$ up to graded isomorphism, analogous considerations ensure that our concept of a coarse moduli space coincides with the original one.  In this situation, graded isomorphism takes on the role of the equivalence relation $\sim$.
\smallskip

\noindent \noindent {\rm \bf{(iii)}}  Clearly, any fine moduli space for $\C$ is a coarse moduli space.   In particular, by Proposition 3.2, either type of moduli space for our problem is an orbit space based on our choice of parametrizing variety (the subvariety $\C(\modlad)$ of $\modlad$ or $\C(\grassbd)$ of $\grassbd$ corresponding to $\C$) modulo the appropriate group action.  From the definition of an orbit space, we thus glean that classification by a coarse moduli space $X$ also yields a one-to one correspondence between the points of $X$ and the isomorphism classes of modules from $\C$.  Concerning fine classification, we moreover observe: If $X$ is a fine moduli space for $\C$, then each isomorphism class from $\C$ is represented by precisely one fibre of the corresponding universal family parametrized by $X$.
\enddefinition   

In essence, the role of a fine or coarse moduli space thus is to not only record parameters pinning down normal forms for the objects in the class $\C$ under discussion, but to do so in an optimally interactive format.  Consequently, under the present angle, the ``effectiveness" of a normal form is measured by the level of universality it carries.  Let us subject some familiar instances to this quality test,
recruiting schoolbook knowledge from the representation theory of  the polynomial algebra $K[t]$. 

\definition{First examples 4.3} 
\smallskip

{\bf (1)}  It is not difficult to check that the family presented in 4.1 is universal for the class $\C$ of non-semisimple modules with dimension vector $(1,1)$ over the Kronecker algebra.  This fact will be re-encountered as a special case of Corollaries 5.2 and 6.7 below. 
\medskip

{\bf (2)} (cf\. \cite{\New, Chapter 2}) Let $d$ be an integer $\ge 2$.  Suppose that $V$ is a $d$-dimensional $K$-space, and $\C$ a class of endomorphisms of $V$.  In other words, we are considering a class of $d$-dimensional modules over $K[t]$.  Rephrasing the above definition of a family of modules, we obtain:  A family from this class, parametrized by a variety $X$ say, is a vector bundle of rank $d$, together with a bundle endomorphism $\delta(t)$ that induces endomorphisms from the class $\C$ on the fibers.  The equivalence relation to be considered is similarity in the usual sense of linear algebra.  

An immediate question arises:  Does the full class $\C$ of endomorphisms of $V$ have a coarse or fine moduli space?

Given that our base field $K$ is algebraically closed, we have Jordan normal forms which are in one-to-one correspondence with the similarity classes.  So the first question becomes:  Can the invariants that pin down the normal forms be assembled to an algebraic variety?  The fact that the block sizes in JNFs are positive integers  --  that is, are discrete invariants  --  while $M_{d \times d}(K)$ is an irreducible variety, does not bode well.  Indeed, one readily finds that all conjugacy classes in $\End_K(V)$ encoding non-diagonalizable endomorphisms (= non-semismple $K[t]$-modules) fail to be closed;  indeed, the Zariski-closure of any such class contains the diagonalizable endomorphism with the same eigenvalues and multiplicities.  Consequently, the full collection of endomorphisms of $V$ does not even have a coarse moduli space.  

 If one restricts to the class $\C$ of diagonalizable endomorphisms of $V$, there is a coarse moduli space; this orbit space $\C /\!/ \GL(V)$ is isomorphic to $K^d$ and records the coefficients of the characteristic polynomial (disregarding the leading coefficient).  But there is no universal family for the problem, so the coarse moduli space fails to be fine in this case.  For a proof, see e.g. \cite{\New, Corollary 2.6.1}.  On the other hand, if one further specializes to the cyclic endomorphisms, i.e., 
$$\C = \{ f \in \End_K(V) \mid f \text{\ corresponds to a cyclic\ } K[t]\text{-module}\},$$
 one finally does obtain a fine moduli space, namely $K^d$; a universal family for the endomorphisms in $\C$ traces their rational canonical forms.  
 \medskip

{\bf (3)}  Riemann's celebrated classification of smooth projective curves of fixed genus over $\CC$ is implemented by a {\it coarse\/} moduli space, which fails to be fine.  This appears to be the situation prevalent in sweeping classification results in algebraic geometry. 
\enddefinition 
\smallskip 

To return to the representation theory of a {\it finite\/} dimensional algebra $\la$:
 Two strategies have emerged to draw profit from the concepts of coarse or fine moduli spaces in this context.  In line with the conclusion of Section 3, each of them reduces the focus to suitable subclasses of the full class of modules with fixed dimension vector.  However, they are based on different expectations, and the dichotomy is paralleled by different techniques.  In the following, we sketch both of these methods and provide sample results.

\head 5.  Approach A:  King's adaptation of Mumford stability:  Focusing on the objects which are (semi-)stable relative to a weight function  \endhead

As the caption indicates, this approach builds on the affine Parametrization A of Section 2. Given that there are already two survey articles recording it, by Geiss \cite{\Geiss} and Reineke \cite{\Reinsurv}, we will be comparatively brief and refer to the existing overviews for technical detail and further applications. 

The strategy under discussion was originally developed for the purpose of classifying certain geometric objects (vector bundles over certain projective varieties, in particular)  subject to the following, a priori unfavorable, starting conditions:  The equivalence classes of the objects are in bijective correspondence with the orbits of a reductive group action on an affine parametrizing variety, but closed orbits are in short supply.  This is precisely the obstacle we encountered at the end of Section 3 relative to $\modlad$ with its $\GL(\bd)$-action.  As a consequence, the attempt to construct an orbit space from invariant regular functions on the considered variety, on the model of Theorem 3.1, is doomed.  The idea now is to use more regular functions, rather than just the classical invariants (constant on the orbits), loosening their tie to the group action to a controllable extent:  namely, to use all regular functions which are  {\it semi-invariant relative to a character of the acting group\/}.  In tandem, one pares down the original variety to an action-stable subvariety with a richer supply of (relatively) closed orbits.   In a nutshell:  One allows for a larger supply of regular functions to palpate a curtailed collection of orbits.  We follow with a somewhat more concrete outline. First we sketch the original GIT-scenario without including the general definitions of (semi)stability and $S$-equivalence.  Then we specialize to the variety $\modlad$ with its $\GL(\bd)$-action and fill in the conceptual blanks, using King's equivalent characterizations of stability and semistability for this case.  (For more precision on the general case, see also \cite{\Der}.)

The typical scenario to which this strategy applies is as follows:  Namely, a finite dimensional $K$-vectorspace $V$ (for example, $V = \modlad$, where $\la = KQ$ is a hereditary algebra), together with a reductive algebraic group $G$ which operates linearly on $V$.  Then a regular function $V \rightarrow K$ is called a {\it semi-invariant\/} for the action in case there exists a character $\chi:  G \rightarrow K^*$ such that $f(g.x) = \chi(g) f(x)$ for all $g \in G$ and $x \in V$.  Next, one singles out a subvariety $V^{\text{st}}$ of $V$ whose $G$-orbits are separated by $\chi$-semi-invariants;  the points of $V^{\text{st}}$ are called ``$\chi$-stable".  In addition, one considers a larger subvariety $V^{\text{sst}}$ whose points are separated by semi-invariants modulo a somewhat coarser (but often still useful) equivalence relation, labeled $S$-equivalence (``$S$" for ``Seshadri"); the points of $V^{\text{sst}}$ are dubbed ``$\chi$-semistable".  More accurately,  the $S$-equivalence classes of $\chi$-semistable points are separated by semi-invariants of the form $\chi^m$ for some $m \ge 0$. The motivation for this setup lies in the following consequences:  The collection $V^{\text{sst}}$ of semistable points is an open (possibly empty) subvariety of $V$ which allows for a categorical quotient that classifies the orbits in $V^{\text{sst}}$ up to $S$-equivalence.  The subset $V^{\text{st}}$ of stable points in $V$ is in turn open in $V$ and far better behaved from our present viewpoint:  It frequently permits even a geometric quotient modulo $G$.  As is to be expected, the quotient of $V^{\text{sst}}$ modulo the $G$-action is constructed from semi-invariants, namely as Proj of the following graded ring of semi-invariant functions:  $\bigoplus_{m \ge 0} K[V]^{\chi^m}$, where $K[V]^{\chi^m}$ is the $K$-subspace of the coordinate ring $K[V]$ consisting of the polynomial functions which are semi-invariant relative to $\chi^m$. 
 
For the module-theoretic scenario that resembles the GIT-template the most closely, King's adaptation of the outlined strategy has been  the most successful.  It is the case of a hereditary algebra $\la = KQ$.  In this situation, $V = \modlad$ is a finite-dimensional $K$-vector space, and the reductive group action is the $G = \GL(\bd)$-conjugation action.  In particular, King showed that $\chi$-(semi)stability, for a character $\chi$ of $\GL(\bd)$, translates into a manageable condition for the modules represented by the $\chi$-(semi)stable points; see below.  He then proceeded to carry over the technique to arbitrary finite dimensional algebras $\la = KQ/I$.

As is well-known, the characters of $\GL(\bd)$ are in natural correspondence with the maps $Q_0 \rightarrow \ZZ$ (see \cite{\Der}, for instance).  Namely, every character $\chi$ is of the form $\chi(g) = \prod_{i \in Q_0} \det(g_i)^{\theta(i)}$ for a suitable map $\theta: Q_0 \rightarrow \ZZ$;  conversely, all maps of this ilk are obviously characters.  Starting with the additive extension $\ZZ^{Q_0} \rightarrow \ZZ$  of such a map  --  called by the same name  --   one lets $\chi_{\theta}$ be the corresponding character of $\GL(\bd)$.  As was proved by King \cite{\King, Theorem 4.1}, a point in $\modlad$ is $\chi_{\theta}$-semistable in the GIT-sense if and only if the corresponding module $M$ satisfies $\theta(\underline{\dim}\, M) = 0$ and $\theta(\underline{\dim}\, M' )\ge 0$ for all submodules $M'$ of $M$; stability requires that $\underline{\dim}\, M$ belong to the kernel of $\theta$ and $\theta(\underline{\dim}\, M') > 0$ for all proper nonzero submodules $M'$ of $M$. For convenience, one also refers to a module $M$ as $\theta$-(semi)stable if it is represented by a $\theta$-(semi)stable point in $\modlad$.  (Note: The function $\theta$ is called a {\it weight\/} by Derksen, a {\it stability\/} by Reineke.) 

Since the sets of $\theta$-semistable, resp\. of $\theta$-stable, points in $\modlad$ are open in $\modlad$, the classes $\C$ of $\theta$-semistable, resp\. $\theta$-stable modules satisfy the blanket hypothesis we imposed in our definitions of a fine or coarse moduli space for $\C$. 

 This setup yields the following:
\smallskip

\proclaim{Theorem 5.1} {\rm {(}}{\rm see} \cite{\King, Propositions 4.3, 5.2 and 5.3}{\rm {)}}.  Let $\la = KQ/I$.
\smallskip
$\bullet$  The $\theta$-semistable objects in $\lamod$ form a {\rm {(}}full\/{\rm {)}} abelian subcategory of $\lamod$ in which all objects have Jordan-H\"older series.  The simple objects in this category are precisely the $\theta$-stable modules.  Two semistable objects are $S$-equivalent precisely when they have the same stable composition factors. 
\smallskip

$\bullet$ The $\theta$-semistable modules of a fixed dimension vector $\bd = (d_1, \dots, d_n)$ have a coarse moduli space $\M^{\text{sst}}_\la (\bd, \theta)$  which classifies them up to $S$-equivalence.  This coarse moduli space is projective and contains, as an open subvariety, a coarse moduli space $\M^{\text{st}}_\la (\bd, \theta)$ classifying the $\theta$-stable modules up to isomorphism.  $\M^{\text{st}}_\la (\bd, \theta)$ is a fine moduli space provided that \, gcd$(d_1, \dots, d_n) = 1$. 
\endproclaim 

Evidently, conditions guaranteeing that $\M^{\text{sst}}_\la (\bd, \theta)$ $\bigl($resp., $\M^{\text{st}}_\la (\bd, \theta) \bigr)$ be nonempty are among the most pressing points to be addressed.  We pinpoint one of the lucky situations, where an effective weight function $\theta$ is easy to come by.  It concerns the classification of local modules (= modules with a simple top), when the quiver of $\la$ has no oriented cycles.  

\proclaim{Corollary 5.2}{\rm  [Crawley-Boevey, oral communication]}  Suppose $\la = KQ/I$, where $Q$ is a quiver without oriented cycles, $T$ is a simple $\la$-module and $\bd$ a dimension vector.  Then the {\rm{(}}local\/{\rm{)}} modules with top $T$ and dimension vector $\bd$ have a fine moduli space, $\M^{\text{st}}_\la (\bd, \theta)$ for a suitable weight function $\theta$, which classifies them up to isomorphism.
\endproclaim

\demo{Proof}  Suppose that $T = S_1$.   Let $\theta: Q_0 \rightarrow \ZZ$ be defined by $\theta(e_j) =  1$ for $j > 1$ and $\theta(e_1) = - \sum_{2 \le j \le n} d_j$.  Then clearly  the modules addressed by the corollary are the $\theta$-stable ones, and King's theorem applies. \qed \enddemo    
\medskip

\centerline {{\bf Pros and cons of Approach A}}
\smallskip
\noindent  {\bf Pros:}

$\bullet$  This tactic always leads to a moduli space if one extends the notion to the empty set.  Indeed, for any choice of weight function, existence of coarse, resp\. fine moduli spaces, for the corresponding semistable, resp\. stable modules, is guaranteed by GIT.
\smallskip

$\bullet$  Since this method has proved very effective for vector bundles on non-singular projective curves, a large arsenal of methods for analyzing the resulting moduli spaces has been developed.  This includes cohomology groups and their Betti numbers, as well as cell decompositions.  (Interesting adaptations to the representation-theoretic setting of techniques developed towards the understanding of vector bundle moduli can be found in the work of Reineke, e.g., in \cite{\Rein}.) 
\smallskip

$\bullet$  The spotlight placed on semi-invariant functions on $\modlad$ by this method appears to have reinforced research into rings of semi-invariants, a subject of great interest in its own right.  
\medskip

\noindent {\bf Cons:}

$\bullet$  How to judiciously choose weight functions is a tough problem.  In this context, a weight  function $\theta: Q_0 \rightarrow \ZZ$ merits the attribute ``good" if one is able to secure a rich supply of $\theta$-stable representations, next to a solid grasp of ``who they are".  ($\M^{\text{st}}_\la (\bd, \theta)$ may be empty.)  There are not (yet) any systematic responses to this problem, beyond some partial insights in the hereditary case. 

\smallskip  

$\bullet$  In general, the $\theta$-(semi)stable modules do not have descriptions in structural terms that turn them into representation-theoretically distinguished classes.  
\smallskip

$\bullet$  The stable modules typically have large orbits, which means that the moduli space $\M^{\text{st}}_\la (\bd, \theta)$ is unlikely to capture boundary phenomena in the geometry of $\modlad$. 
\smallskip

$\bullet$  This refers to a weight function $\theta$ such that $\M^{\text{sst}}_\la (\bd, \theta)$ is nonempty:  The fact that it is typically difficult to interpret $S$-equivalence in representation-theoretic terms detracts  --  at least for the moment  --   from the value of the existence of coarse moduli spaces that classify the semistable modules up to  this equivalence.  
\bigskip

\noindent{\it Exploring and addressing these problems:}   Here is a selection of insights for the special case where $\la = KQ$:
\smallskip

$\bullet$ Existence of a weight function $\theta$ with the property that $\M^{\text{st}}_\la (\bd, \theta) \ne \varnothing$ is equivalent to $\bd$ being a Schur root (see \cite{\King, Proposition 4.4}).  In fact, stability of a module $M$ relative to some weight function forces $M$ to be a Schurian representation, that is, to have endomorphism ring $K$.  Since the Schurian representations with dimension vector $\bd$ clearly have maximal orbit dimension in $\modlad$, the union of the $\M^{\text{st}}_\la (\bd, \theta)$, where $\theta$ traces different weight functions, is contained the open sheet of $\modlad$ (for sheets, see, e.g., \cite{\Kra}).  The fact that the variety $\M^{\text{st}}_\la (\bd, \theta)$ is always smooth (see \cite{\King, Remark 5.4}) once more points to absence of boundary phenomena.  
\smallskip

$\bullet$  Given a Schur root $\bd$, there is in general no choice of $\theta$ such that all Schur representations $M$ of $\la$ with dimension vector $\bd$ are $\theta$-stable.  In fact, for a given Schurian representation, there need not be any weight function $\theta$ making it $\theta$-stable (see \cite{\Rein, Section 5.2}, where the $5$-arrow Kronecker quiver is used to demonstrate this). 
\smallskip

$\bullet$   On the positive side: Given $Q$, $\theta$ and $\bd$, Schofield's algorithm in \cite{\Scho} permits to decide whether $\bd$ is a Schur root of $KQ$ and, if so, whether there is a $\theta$-stable Schurian representation of dimension vector $\bd$; cf\. \cite{\King, Remarks 4.5, 4.6}.  Furthermore, Reineke developed a recursive procedure for deciding whether $\M^{\text{sst}}_\la (\bd, \theta) \ne \varnothing$;  see \cite{\Rein, Section 5.3} for an outline. (The argument is based on an adaptation to the quiver scenario of results due to Harder and Narasimhan and provides a specific instance of one of the plusses listed above.)   See also the work by Adriaenssens and Le Bruyn \cite{\Adrialebruyn} on assessing the supply of $\theta$-(semi)stable modules with dimension vector $\bd$.      
\bigskip

{\it Pointers for further reading\/}:  Over the past 20 years, this angle on moduli of representations has inspired an enormous amount of research, with interesting results not only directly targeting moduli spaces, but also rings of semi-invariants of the varieties $\modlad$ in their own right, next to surprising applications, for instance to Horn's Problem.  Let us just mention a (necessarily incomplete) list of further contributors:  Chindris, Crawley-Boevey, de la Pe\~na, Derksen, Geiss, Hille, Le Bruyn, Nakajima, Procesi, Reineke, Schofield, Van den Bergh, Weyman.

\head  6.  Approach B.  Slicing $\lamod$ into strata with fixed top  \endhead

In the following, we rely on the projective parametrization introduced in Section 2.

Instead of using stability functions to single out classifiable subvarieties of $\modlad$, we start by partitioning $\grassbd$ into finitely many locally closed subvarieties, based on module-theoretic invariants.  The primary slicing is in terms of tops.  Let $T \in \lamod$ be semisimple.  The restriction to modules with fixed top $T$ has an immediate payoff.  Namely, the locally closed subvariety 
$$\biggrass^T_{\bd} = \{C \in \grassbd \mid \top(\boldP / C) = T\}$$ 
of $\grassbd$ may be replaced by a projective parametrizing variety, $\grasstbd$, which has far smaller dimension in all interesting cases.  In fact, the pared-down variety $\grasstbd$ appears to go part of the way towards a quotient of $\GRASS^T_{\bd}$ by its $\aut_\la(\boldP)$-action.  In many instances, moduli spaces for substantial classes of representations with fixed top $T$ will, in fact, be identified as suitable closed subvarieties of $\grasstbd$.  This is for instance true in the local case addressed in Corollary 5.2:  The fine moduli space $\M^{\text{st}}_\la (\bd, \theta)$, guaranteed by Approach A in that case, equals $\grasstbd$; see Corollary 6.7 below for justification. 
\smallskip

Following the tenet ``smaller is better", we fix a projective cover $P$ of $T$, to replace the projective cover $\boldP$ of $\bigoplus_{1 \le i \le n} S_i^{d_i}$. Since we are restricting our focus to modules with top $T$, this projective cover suffices.  Accordingly, we consider the subset 
$$\grasstbd = \{ C \in \Gr(\dim P - |\bd|\, , P) \mid \,  _\la C \text{\ is a submodule of\ }  _\la JP \text{\ and\ } \underline{\dim}\, P/C = \bd \}$$
of the classical Grassmannian consisting of the ($\dim P - |\bd|$)-dimen\-sional $K$-subspaces of the $K$-vector space $JP$. 
Clearly, $\grasstbd$ is in turn a closed subvariety of the subspace Grassmannian $\Gr(\dim P - |\bd|\, , P)$, and hence is projective.  Moreover, the natural action of the automorphism group  $\autlap$ on $\grasstbd$ once more provides us with a one-to-one correspondence between the set of orbits on one hand and the isomorphism classes of $\la$-modules with top $T$ and dimension vector $\bd$ on the other.  Evidently, we have the same semi-direct product decomposition of the acting group as before:  $\autlap \cong \aut_\la(T) \ltimes \unirad$, where $\unirad = \{\id_P + h \mid h \in \hom_\la(P, JP)\}$ is the unipotent radical of $\autlap$. 

The main reason for expectations of a gain from this downsizing is as follows:  The semi-direct product decomposition of the acting automorphism group, in both the big and small scenarios, invites us to subdivide the study of orbit closures into two parts. It does, indeed, turn out to be helpful to separately focus on orbits under the actions of the semidirect factors, and it is foremost the size of the reductive factor group, $\aut_\la (T) = \aut_\la (P/JP)$ resp\. $\aut_\la(\boldP  / J\boldP)$, which determines the complexity of this task.  (In Section 7, it will become apparent why the action of the unipotent radical is easier to analyze.)   As a consequence, it is advantageous to pass from the big automorphism group $\aut_\la(\boldP)$ to one with reductive factor group as small as possible.   Corollaries 6.6, 6.7 and Proposition 7.2, in particular, attest to the benefits that come with a simple, or at least squarefree, top.

\definition{6.1. Preliminary examples}  {\bf (1)}  Let $\la = KQ$, where $Q$ is the generalized Kronecker quiver with $m \ge 2$ arrows from a vertex $e_1$ to a vertex $e_2$.  Moreover, choose $T = S_1$ and $\bd = (1,1)$.  Then $P = \la e_1$,  $\grasstbd \cong \PP^{m - 1}$, and the $\autlap$-orbits are singletons.  Thus $\grasstbd$ is an orbit space.  Corollary  5.2 guarantees a fine moduli space classifying the modules with top $T$ and dimension vector $\bd$ up to isomorphism, and hence $\grasstbd$ coincides with this moduli space.
\smallskip

\noindent {\bf (2)} Next, let $\la = KQ/I$, where $Q$ is the quiver
$$\xymatrixrowsep{3pc}\xymatrixcolsep{6pc}
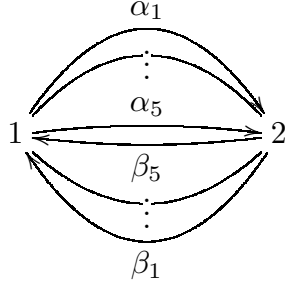
\xymatrix{
1 \ar@/^2.8pc/[r]^{\alpha_1} \ar@{}@/^2.1pc/[r]|{\vdots} \ar@/^0.25pc/[r]^{\alpha_5}  &2  \ar@/^2.8pc/[l]^{\beta_1} \ar@{}@/^1.8pc/[l]|{\vdots} \ar@/^0.25pc/[l]^{\beta_5}
}$$
\noindent and $I$ is the ideal generated by the $\beta_i \alpha_j$ for $i \ne j$ and all paths of length $3$. Again choose $T = S_1$.  For $\bd = (d'_1 +1, d_2)$ with $d'_1, d_2 \le 5$, we obtain the following distinct outcomes concerning $\grasstbd$:
If $d'_1 > d_2$, then $\grasstbd$ is empty. If $d'_1 = d_2$, then $\grasstbd  \cong \Gr(5 - d'_1, K^5)$ $\cong$ $\Gr(d'_1, K^5)$. If $d'_1 < d_2$, then $\grasstbd \cong \flag(5 - d'_1, 5 - d_2, K^5)$, where the latter denotes the variety of partial subspace flags $K^5 \supseteq U_1 \supseteq U_2 $ with $\dim  U_1 =  5 - d'_1$ and $\dim U_2 = 5 - d_2$.  As a consequence of Corollary 6.7 below, we will find that, in either case, $\grasstbd$ is a fine moduli space classifying the modules with top $T$ and dimension vector $\bd$ up to isomorphism. 
\smallskip

\noindent {\bf (3)}  Finally, let $\la = KQ/I$, where $Q$ is the quiver $\xymatrixcolsep{2.5pc} \xymatrix{ 1 \ar@(ul,dl)_{\alpha} \ar[r]^{\beta} &2 }$, and $I = \langle \alpha^2 \rangle$.  For $T = S_1$ (hence $P = \la e_1$) and $\bd = (2,1)$, we obtain $\grasstbd \cong \PP^1$.  From Section 3, we glean that the modules with top $T$ and dimension vector $\bd$ do not even have a coarse moduli space classifying them up to isomorphism.   Indeed, the $\autlap$-orbit of the point $C = \la \beta \in \grasstbd$ is a copy of $\AA^1$, and consequently fails to be closed in $\grasstbd$.   On the other hand, the modules in $\grasstbd$ are classifiable in naive terms --  up to isomorphism, there are only two of them after all.  In order to obtain the benefits of a fine classification in the strict sense, however, one needs to stratify $\grasstbd$ further into segments with fixed radical layerings. In the present example, this is a trivial stratification into $\AA^1$ and a singleton.
\enddefinition

In the present smaller setting, the transfer of information between the projective and the affine parametrizing varieties follows the same pattern as in the big (described in Proposition 3.2).  Clearly, the counterpart of $\grasstbd$ in the affine setting is the subvariety $\toptbd$ of $\modlad$ which consists of the points that represent modules with top $T$.  
Observe that $\grasstbd$ records the same geometric information as $\biggrass^T_{\bd}$, just in a less redundant format.

\proclaim{Proposition 6.2. Information transfer revisited}
Let $\Psi$ be the bijection  
$$\{ \autlap\text{-stable subsets of\ } \grasstbd \} \rightarrow
\{\GL(\bd)\text{-stable subsets of\ } \toptbd \}$$
extended from the one-to-one correspondence between sets of orbits
which assigns to an orbit $\autlap.C$ of $\grasstbd$ the $\GL(\bd)$-orbit
of $\toptbd$ that represents the isomorphism class of $P/C$.  Once again, $\Psi$ is an inclusion-preserving bijection which preserves and reflects
openness, closures, connectedness, irreducibility, and types of
singularities.  Moreover, it preserves categorical and geometric quotients of $\autlap$-stable subvarieties of $\grasstbd$, as well as orbit spaces for the $\autlap$-action.  The inverse $\Psi^{-1}$ has analogous preservation properties.
\endproclaim

Next we present a selection of results addressing existence and, if pertinent, properties of fine or coarse moduli spaces for:  (I) The modules that do not admit any proper top-stable degenerations.  (II)  The graded modules with fixed top and dimension vector over an algebra $\la = KQ/I$, where $I$ is a homogeneous ideal.

\bigskip

\subhead  (I)  The modules which are degeneration-maximal among those with fixed top  \endsubhead 
\smallskip

What are they?

Let $M$ and $M'$ be $\la$-modules with dimension vector $\bd$.  Recall that $M'$ is a {\it degeneration\/} of $M$ in case the $\GL(\bd)$-orbit in $\modlad$ that corresponds to $M'$ is contained in the closure of the $\GL(\bd)$-orbit corresponding to $M$.  By Proposition 3.2, this amounts to the same as postulating that the $\aut_\la(\boldP)$-orbit representing $M'$ in $\grassbd$ be contained in the closure of the $\aut_\la(\boldP)$-orbit representing $M$.  We write $M \degen M'$ to communicate this connection between the orbits, and observe that $\degen$ defines a partial order on isomorphism classes of modules.  Intuitively, one may think of the degenerations of $M$ as a collection of modules that document a successive unraveling of the structure of $M$, following a geometry-guided instruction set;  this viewpoint is buttressed by examples.

Note:  Whereas in the present context  --  the pursuit of moduli spaces  --  the typically enormous sizes of orbit closures in module varieties is a priori an obstacle, a shifted viewpoint makes a virtue out of necessity.  One way of organizing the category $\lamod$ is to break it up into posets of (isomorphism classes of) degenerations of individual modules, and to analyze these posets in their own right; this direction has, in fact, moved to the mainstream of research.  Along a related line, it is profitable to take aim at those modules in a specified subvariety $X$ of $\modlad$ which are distinguished by having the same ``height" (or ``depth") relative to the degeneration order within $X$.  (Observe that, for given $\bd$, the lengths of chains of degenerations of modules with dimension vector $\bd$ are bounded from above by $|\bd| - 1$; we follow the Romans and start with $0$ in counting chain lengths.)  This is, in fact, the tack we are taking in this subsection.  For background on the extensive theory of degenerations we refer the reader to work of Bobinski, Bongartz, Riedtmann, Schofield, Skowronski, Zwara, and the author, for instance.  Three seminal articles provide a good point of departure: \cite{\Riedt}, \cite{\Bongadv}, \cite{\Zwara}.  

The representations which are maximal under $\degen$ in $\lamod$ do not hold much interest.  It is easy to see that, given any submodule $U$ of a module $M$, the direct sum $U \oplus M/U$ is a degeneration of $M$.  Hence, for any dimension vector $\bd$, there is, up to isomorphism, exactly one module which is degeneration-maximal among the modules with that dimension vector, namely the semisimple module $\bigoplus_{1 \le i \le n} S_i^{d_i}$.  By contrast, there is usually a plethora of degenerations of $M$ which are maximal among the degenerations that have the same top as $M$  (see 6.4 below).  On the other hand, Theorem 6.3 below guarantees that they nonetheless always have a fine moduli space classifying them up to isomorphism.   

Since we are focusing on modules with fixed top $T$, it is the orbit closure of a module $M$ in $\grasstbd$ (resp., $\toptbd$) that is relevant for the moment.  Accordingly, we refer to $M'$ as a {\it top-stable degeneration} of $M$ in case $M \degen M'$ and $\top(M) = \top(M')$.  Clearly, $M = P/C$ with $C \in \grasstbd$ is degeneration-maximal among the modules with top $T$ (meaning that $M$ has no proper top-stable degeneration) precisely when the $\autlap$-orbit of $C$ is closed in $\grasstbd$. 

\proclaim{Theorem 6.3} {\rm(see \cite{\DHZW, Theorem 4.4 and Corollary 4.5})}  For any semisimple $T \in \lamod$, the modules of dimension vector $\bd$ which are degeneration-maximal among those with top $T$ have a fine moduli space, $\maxmoduli$, that classifies them up to isomorphism.

The moduli space $\maxmoduli$ is a closed subvariety of $\grasstbd$, and hence is projective.  

In particular, given any module $M$ with dimension vector $\bd$ whose top is contained in $T$, the closed subvariety of $\maxmoduli$ consisting of the points that correspond to degenerations of $M$ is a fine moduli space for the maximal top-$T$ degenerations of $M$.
\endproclaim

Observe that $\top(M) \subseteq \top(M')$ whenever $M \degen M'$.  By the theorem, we hit new classifiable strata in the hierarchy of degenerations of $M$ as we successively enlarge the allowable top.

The moduli space $\maxmoduli$ is located in $\grasstbd$ as follows:  First one zeroes in on the subvariety $\M$ of $\grasstbd$ consisting of the closed orbits (that is, on the orbits of the target class of modules).  On $\M$, the $\unirad$-action is trivial, but $\aut_\la(T)$ will still operate with orbits of arbitrarily high dimension in general.  However, if we pick a Borel subgroup $H$ of $\aut_\la(T)$ and cut $\M$ back to the closed subvariety of all points that have a stabilizer containing $H$, we arrive at an incarnation of $\maxmoduli$. 

The following concomitant result provides evidence for the representation-theoretic richness of the classes of representations addressed by Theorem 6.3.  The construction used has predecessors in \cite{\GeomI, Theorem 6} and \cite{\Hille, Example}.

\proclaim{Satellite result 6.4} {\rm {(}}{\rm see} \cite{\DHZW, Example 5.4}{\rm {)}} Every projective variety is isomorphic to $\maxmoduli$ for some choice of $\la$, $T$, and $\bd$. \endproclaim

A crucial ingredient of the proof of Theorem 6.3 consists of normal forms of the modules without proper top-stable degenerations.  In fact, the shape of their normalized projective presentations is both of independent interest and guides the explicit construction of universal families.   The reformulation of absence of proper top-stable degenerations under $(1)$ below is due to projectivity of the variety $\grasstbd$.  (By definition, a closed subgroup $H$ of a linear algebraic group $G$ is parabolic precisely when the geometric quotient $G/H$ is a projective variety.)

\proclaim{Theorem 6.5}  {\rm {(}}{\rm see} \cite{\DHZW, Theorem 3.5}{\rm {)}} Let $M$ be a module with  dimension vector $\bd$ and top $T = \bigoplus_{1 \le i \le n} S_i^{t_i}$.  Moreover, let $C$ be a point in $\boldgrasstd$ such that $M \cong P/C$.  
Then the following statements are equivalent:  

\roster
\item"{\rm (1)}" $M$ has no proper top-stable degenerations, i.e., the stabilizer subgroup $\stab_{\autlap} (C)$ is a parabolic subgroup of $\autlap$.
\item"{\rm (2)}"  $M$ satisfies these two conditions:
\itemitem{$\bullet$}  $M$ is a direct sum of local modules, say $M =
\bigoplus_{1
\le i \le n} \bigoplus_{1 \le j \le t_i} M_{ij}$, where $M_{ij} \cong \la
e_i/C_{ij}$ with the following additional property:
For each $i \le n$, the $C_{ij}$ are linearly ordered
under inclusion.
\itemitem{$\bullet$} $\dim_K \hom_\la(P, JM) = \dim_K \hom_\la(M,JM)$.
\endroster

\noindent If conditions {\rm (1) -- (3)} are satisfied, then $\unirad$ stabilizes $C$, and $\autlap.C = \autlat.C$ is isomorphic to a direct product of partial flag
varieties $\F_i$,  where $\F_i$ depends only on the number of distinct left ideals in the family $(C_{ij})_{j \le t_i}$ and their multiplicities.
\endproclaim

The dimension condition in statement (2) of Theorem 6.5 has the following interpretation:  It means that the first syzygy of $M$ is invariant under all homomorphisms $P \rightarrow JP$. 

The upcoming corollaries rest on the following combination of Theorem 6.3 with Section 3.  It shows that one of the cons we listed in connection with Method A arises in Method B as well:  Namely, for large tops $T$, few {\it closed\/} subvarieties $X$ of $\toptbd$ correspond to classes of modules permitting a fine moduli classification.   In other words, ``most" classifications of this ilk are expected to target only generic classes of modules, thus circumventing ``boundary phenomena". 
\smallskip

\noindent{\bf Consequence concerning the classifiability of closed subsets of $\toptbd$.}  Let $X$ be a closed $\autlap$-stable subvariety of $\grasstbd$ and $\C$ the class of modules represented by the orbits of $X$.  Then there is a fine (equivalently, a coarse) moduli space classifying the modules in $\C$ up to isomorphism if and only if $\C$ consists of modules that are degeneration-maximal among those with top $T$.   
\smallskip

In Corollaries 6.6 and 6.7, it does not affect the outcome of the moduli problem whether we fix a dimension vector or else fix only the total dimension of the modules considered.  We opt for the latter, since this leads to smoother statements. To that end, we slightly upgrade our notation:  For any positive integer $d$, we denote by $\grasstd$ the union of the varieties $\grasstbd$ where $\bd$ ranges over the dimension vectors with $|\bd | = d$; that is, $\grasstd =  \{ C \in \Gr(\dim P - d, P) \mid  \, _\la C \text{\ is a submodule of\ }  _\la JP \}$.  

\proclaim{Corollary 6.6}   {\rm {(}}{\rm see} \cite{\classifying}{\rm {)}}  Suppose that $T$ is a squarefree semisimple module and $d \in \NN$.  Then the following statements are equivalent:
\roster
\item"{\bf (a)}"   The modules with top $T$ and dimension $d$ have a fine moduli space classifying them up to isomorphism.       
\item"{\bf (b)}"   The modules with top $T$ and dimension $d$ have a coarse moduli space classifying them up to isomorphism.
\item"{\bf (c)}"  The submodules of $JP$ of dimension $\dim P - d$ are invariant under all endomorphisms of $P$.
\item"{\bf (d)}"   $\grasstd$ is a fine moduli space classifying the $d$-dimensional modules with top $T$ up to isomorphism.
\endroster
\endproclaim

\demo{Deducing the corollary from Theorem 6.5}  The implications (d)$\implies$(a)$\implies$(b) are clear.  Re ``(b)$\implies$(c)":  From (b) we infer that, for each $C \in \grasstd$, the module $M = P/C$ is without proper top-stable degenerations.  Hence, Theorem 6.5 yields the invariance of $C$ under homomorphisms $P \rightarrow JP$, as noted above.  For ``(c)$\implies$(d)", observe that (c) forces all $\autlap$-orbits of $\grasstd$ to be singletons and thus makes all modules $P/C$ with $C \in \grasstd$ degeneration-maximal among the modules with top $T$.  Hence $\grasstd$ is an orbit space whose fibers are singletons.  One may now either invoke Theorem 6.3 or verify that the tautological family parametrized by $\grasstd$ is universal. \qed \enddemo    

We add an offshoot to the previous corollary.  It is a mild extension of Corollary 5.2.

\proclaim{Corollary 6.7}  Suppose that $T$ is simple and that the only occurrences of $T$ in $JP$ are in the socle.  Then the equivalent conditions above are satisfied for all $d$.  Consequently, each $\grasstd$ is a fine moduli space classifying the $d$-dimensional modules with top $T$ up to isomorphism. 
\endproclaim

Let us return to Example (2) in 6.1.  That $\grasstbd$ is a fine moduli space in this instance as well is a special case of Corollary 6.7.
\bigskip

\subhead  (II)  Aiming at the graded modules with fixed top  \endsubhead 
\smallskip

Suppose $\la = KQ/I$, where $I \subseteq KQ$ is a homogenous ideal relative to the path-length grading of $KQ$.  Whenever we speak of graded (left) $\la$-modules we refer to the path-length grading of $\la$.  It is hardly surprising that the additional rigidity encountered in the category of graded modules with homogeneous homomorphisms  (of degree $0$) promotes classifiabilty.  We explore to what extent.

Let $T \in \lamod$ be a semisimple module endowed with the grading that makes it homogeneous of degree $0$.  It is a matter of course that, in addressing graded representations with fixed top $T$ and dimension $d$, we should replace the parametrizing variety $\grasstd$ by a graded incarnation, that is, by
$$\gradgrasstd = \{C \in \grasstd \mid C \text{\ is a homogeneous submodule of\ } JP\};$$
where $P$ stands for the graded projective cover of $T$. 
In tandem, we replace the acting group $\autlap$ by the subgroup $\gradautlap$ consisting of the homogeneous automorphisms of $P$.  One readily confirms that the natural (morphic) action of $\gradautlap$ on $\gradgrasstd$ places the $\gradautlap$-orbits of $\gradgrasstd$ into a canonical one-to-one correspondence with the graded $\la$-modules with top $T$ and dimension $d$.  This setup yields a significant improvement of Corollary 6.7 in the graded situation.

\proclaim{Theorem 6.8} \cite{\BHTone, Theorem 4.1}  Here ``graded" includes ``generated in degree $0$".  

For any simple module $T$ and $d \in \NN$, the $d$-dimensional graded $\la$-modules with top $T$ possess a fine moduli space classifying them up to graded isomorphism.  This moduli space equals $\gradgrasstd$ and, in particular, is a projective variety.
\endproclaim

It is now clear that the finite direct sums of local graded modules are classifiable by fine moduli spaces in segments, namely after the obvious subdivision according to tops and sequences of dimensions of the local summands with fixed top.  This is as far as this kind of ``global" classification can be pushed in the graded case.  The following result attests to a roadblock.  We still include ``generated in degree $0$" when we refer to graded modules.

\proclaim{Theorem 6.9} \cite{\BHTtwo, Theorem 4.2}  Let $T$ be any semisimple module endowed with the obvious grading and $d \in \NN$.  

  If the graded modules with top $T$ and dimension $d$ have a coarse moduli space classifying them up to graded isomorphism, then they are all direct sums of local modules.
\endproclaim

\bigskip

\subhead  (III)  Three easy pieces  \endsubhead 
\medskip

The fine moduli spaces we encountered under (I) and (II) and the corresponding universal families are accessible to algorithmic computation, to the extent that there is an algorithm for determining the distinguished affine cover of these moduli spaces in terms of polynomial equations; it is induced by the distinguished affine cover $\bigl( \grassS \bigr)_{\S}$ of the ambient $\grasstd$; see Section 7.  The restrictions of the targeted universal family to the charts of this cover can in turn be calculated.

The first two pieces, the easiest, illustrate Corollary 6.7 and Theorem 6.8.

\definition{Example 6.10}  Suppose that $J^2 = 0$, i.e., $\la$ is of the form $KQ/I$ where $I$ is generated by all paths of length $2$ in $Q$.  As usual, $Q_0 = \{e_1, \dots, e_n\}$.   Moreover, let $T \in \lamod$ be simple, say $T = S_1$, and $d \in \NN$.  Then the irreducible components of the fine moduli space $\gradgrasstd = \grasstd$ are direct products of classical Grassmannians $\Gr(u_i, K^{v_i})$ for $i \le n$, where the $v_i$ are the numbers of distinct arrows $e_1 \rightarrow e_i$, respectively, and the $u_i \le v_i$ are subject to the equality $1 + \sum_{1 \le i \le n} (v_i - u_i) \, = \, d$. 

Indeed, $J e_1 \cong \bigoplus_{1 \le i \le n} S_i^{v_i}$.  Hence, every $d$-dimensional $\la$-module $M$ with top $T$ has a first syzygy of the form $C = \bigoplus_{1 \le i \le n} W_i$, where $W_i$ is a subspace of $S_i^{v_i}$ whose dimensions add up to $\dim \la e_1 - d$. 
\enddefinition 

\definition{Example 6.11}  let $Q$ be the quiver 
$$\xymatrixcolsep{4pc}
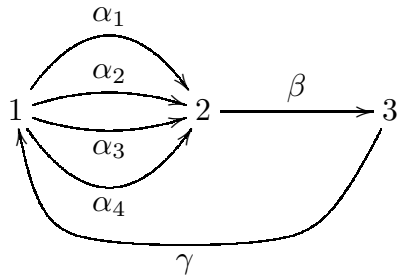
\xymatrix{
1 \ar@/^2pc/[r]^{\alpha_1} \ar@/^/[r]^{\alpha_2} \ar@/_/[r]_{\alpha_3} \ar@/_2pc/[r]_{\alpha_4} &2 \ar[r]^{\beta} &3 \ar@'{@+{[0,0]+(-6,-12)}@+{[0,0]+(-10,-15)}@+{[0,-2]+(6,-15)}@+{[0,-2]+(2,-12)}}[ll]^{\gamma} 
}$$
\noindent and $\la = KQ/I$ where $I \subseteq KQ$ is the ideal generated by  all paths of length $4$.  Moreover, let $T = S_1$ and $\bd = (2, 3, 2)$.  Clearly, all local $\la$-modules are graded, and the fine moduli space classifying the modules with top $T$ and dimension vector $\bd$ up to isomorphism is $\grasstbd \cong \flag(K^4)$.  Indeed, all of the considered modules have radical layering $(S_1, S_2^3, S_3^2, S_1)$ and consequently are of the form $P/C$ with $P = \la e_1$.  The claim can be read off the graph of $P$, which is
$$\xymatrixrowsep{2pc}\xymatrixcolsep{0.5pc}
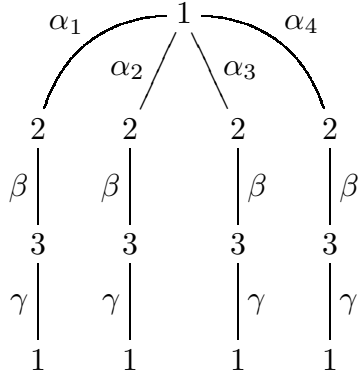
\xymatrix{
&&&1 \edge@/_1pc/[dlll]_{\alpha_1} \edge[dl]_(0.55){\alpha_2} \edge[dr]^(0.55){\alpha_3} \edge@/^1pc/[drrr]^{\alpha_4}  \\
2 \edge[d]_{\beta} &&2 \edge[d]_{\beta} &&2 \edge[d]^{\beta} &&2 \edge[d]^{\beta}  \\
3 \edge[d]_{\gamma} &&3 \edge[d]_{\gamma} &&3 \edge[d]^{\gamma} &&3 \edge[d]^{\gamma}  \\
1 &&1 &&1 &&1
}$$
\noindent (For an informal description of our graphing technique, we refer to \cite{\menace}.) Indeed, we find that $C = W_1 \oplus W_2 \oplus W_3$, where $W_1$ is a $1$-dimensional subspace of $\bigoplus_{i=1}^4 K \alpha_i$, $W_2$ is a $2$-dimensional subspace of $\bigoplus_{i=1}^4 K \beta_i$ which contains $\beta W_1$, and $W_3$ is a $3$-dimensional subspace of $J^3 e_1$ which contains $\gamma W_2$.  
\enddefinition
\smallskip

Next we illustrate Theorem 6.3, focusing on the maximal top-stable degenerations of a single module $M$.  Observe that, in the case addressed below, $M$ is ``close" to being degeneration-maximal among the modules with top $S_1^2$ (indeed, $M$ satisfies all but the last of the conditions in Theorem 6.5(2)). Nonetheless, the poset of top-stable degenerations of $M$ has chains of length $3$.

\definition{Example 6.12}  Let $K = \CC$ and $\la = KQ/I$, where $Q$ is the quiver
$$\xymatrixcolsep{5pc}
 
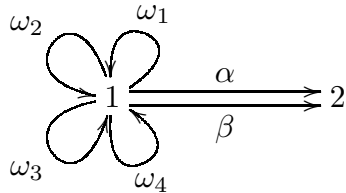
\xymatrix{ 
1 \ar@'{@+{[0,0]+(7,4)}@+{[0,0]+(4,9)}@+{[0,0]+(-1,7)}}_{\omega_1} 
\ar@'{@+{[0,0]+(-2,8)}@+{[0,0]+(-8,8)}@+{[0,0]+(-8,1)}}_{\omega_2} 
\ar@'{@+{[0,0]+(-8,-1)}@+{[0,0]+(-8,-8)}@+{[0,0]+(-2,-8)}}_{\omega_3} 
\ar@'{@+{[0,0]+(-1,-7)}@+{[0,0]+(3.5,-9)}@+{[0,0]+(7,-4)}}_{\omega_4}
\ar[r]<0.5ex>^{\alpha}  \ar[r]<-0.5ex>_{\beta} &2     
 }$$
and $I = \langle \omega_i \omega_j \mid 1 \le i,j \le 4 \rangle + \langle \alpha \omega_i \mid i = 3,4 \rangle + \langle \beta \omega_i \mid i= 1,2 \rangle.$

We take $T = S_1^2$ with projective cover $P = \la z_1 \oplus \la z_2 \cong (\la e_1)^2$, and choose $M = P/ C$ with $C = \la(\alpha + \beta)z_2 + L z_2$, where $L \subseteq \la e_1$ is the left ideal of $\la$ generated by $\alpha \omega_1 + 2 \alpha \omega_2$ and $\beta \omega_3 + 3 \beta \omega_4$; note that $L$ is actually a two-sided ideal.  The module $M$ may be visualized by way of the following ``hypergraph".  For a clean definition of a {\it hypergraph\/} of a module, we refer to \cite{\BHTtwo, Definition 3.9}. 
$$\xymatrixrowsep{2pc}\xymatrixcolsep{1.5pc}
\xymatrix{
 &&&1 \dropup{z_1} \edge@/_2pc/[dlll]_{\alpha} \edge@/_1pc/[dll]_(0.7){\beta} \edge[dl]_(0.65){\omega_1} \edge[d]^(0.6){\omega_2} \edge@/^/[dr]^(0.75){\omega_3} \edge@/^1pc/[drr]^{\omega_4}  && &&&&&1 \dropup{z_2} \edge@/_2pc/[dlll]_{\alpha} \edge@/_/[dlll]^(0.65){\beta} \edge[dl]_(0.65){\omega_1} \edge[d]^(0.6){\omega_2} \edge@/^/[dr]^(0.75){\omega_3} \edge@/^1pc/[drr]^{\omega_4}  \\
2 &2 &1 \edge[d]_{\alpha} &1 \edge[d]_{\alpha} &1 \edge[d]_{\beta} &1 \edge[d]_{\beta} &\bigoplus &2 &&1 \edge[dr]_{\alpha} &1 \edge[d]^{\alpha} &1 \edge[dr]_{\beta} &1 \edge[d]^{\beta}  \\
 &&2 &2 &2 &2 && &&&2 &&2
 }$$
\noindent Take $\bd = \underbardim M = (10,9)$. The moduli space $\maxmoduli$ for the maximal top-T degenerations of $M$ has two irreducible components, both $1$-dimensional.   One is isomorphic to $\C_1 = \{[c_1:c_2:c_3:c_4] \in \PP^3 \mid c_2 = 2 c_1 \text{\ and \ } c_4 = 3 c_3\}$ with the following informally presented universal family: to a point $[c_1: \cdots: c_4] \in \C_1$ we assign the module $P/D$, where 
$$D = \la\bigl(\sum_{j = 1,2} c_j \alpha \omega_j +  \sum_{j = 3,4} c_j \beta \omega_j \bigr) z_1 +   L z_2.$$
  The other irreducible component of $\maxmoduli$ is isomorphic to $\C_2 = \PP^1$; to a point $[c_1:c_2] \in \C_2$ we assign the factor module
$P / D$, where 
$$D = L z_2 + \la(c_1  \alpha \omega_1 z_2 + c_2 \beta \omega_4 z_2).$$ 
We graphically present the generic module for each of these components.
$$\xymatrixrowsep{2pc}\xymatrixcolsep{1.4pc}
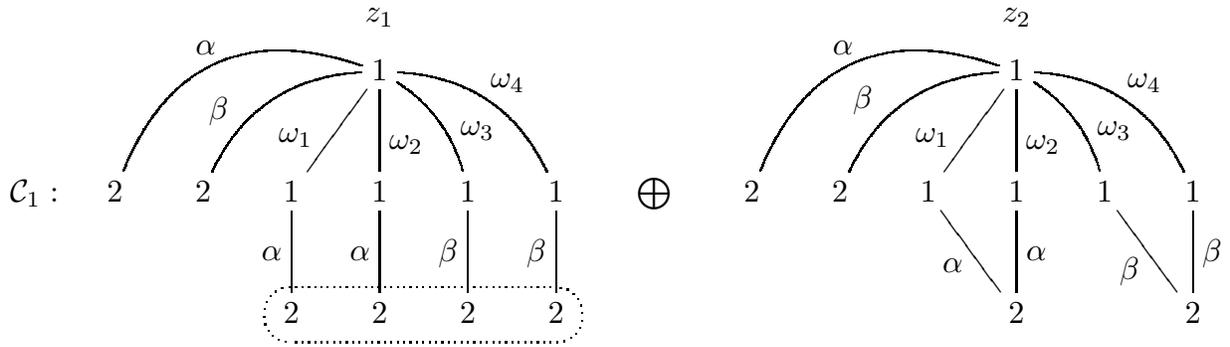
\xymatrix{
 &&&&1 \dropup{z_1} \edge@/_2pc/[dlll]_{\alpha} \edge@/_1pc/[dll]_(0.7){\beta} \edge[dl]_(0.65){\omega_1} \edge[d]^(0.6){\omega_2} \edge@/^/[dr]^(0.65){\omega_3} \edge@/^1pc/[drr]^{\omega_4}  && &&&&&1 \dropup{z_2} \edge@/_2pc/[dlll]_{\alpha} \edge@/_1pc/[dll]_(0.65){\beta} \edge[dl]_(0.65){\omega_1} \edge[d]^(0.6){\omega_2} \edge@/^/[dr]^(0.65){\omega_3} \edge@/^1pc/[drr]^{\omega_4}  \\
\dropvert0{{\Cal C}_1:} &2 &2 &1 \edge[d]_{\alpha} &1 \edge[d]_{\alpha} &1 \edge[d]_{\beta} &1 \edge[d]_{\beta} &\bigoplus &2 &2 &1 \edge[dr]_{\alpha} &1 \edge[d]^{\alpha} &1 \edge[dr]_{\beta} &1 \edge[d]^{\beta}  \\
 &&&2 \levelpool3 &2 &2 &2 & \dropvert{-8}{\left( \sum_{j=1,2} c_j\alpha\omega_j + \sum_{j=3,4} c_j\beta\omega_i \right) z_1 =0 \text{\;\;for\;\;} [c_1:c_2:c_3:c_4] \in {\Cal C}_1} & &&&2 &&2
 }$$ 
\noindent Here the dotted enclosure in the left-hand graph indicates that the corresponding elements are subject to the displayed dependence relation. As we are interested in the generic picture, we assume the $c_j$ to be all nonzero.  The dotted curve enclosing four of the vertices in the bottom row of the preceding graph indicates that the elements $\alpha \omega_1 z_1, \alpha \omega_2 z_1, \beta \omega_3 z_1, \beta \omega_4 z_1$ are linearly dependent  --  the pertinent relation is given below the graph  --  while any three of these elements are linearly independent.  

$$\xymatrixrowsep{2pc}\xymatrixcolsep{1.45pc}
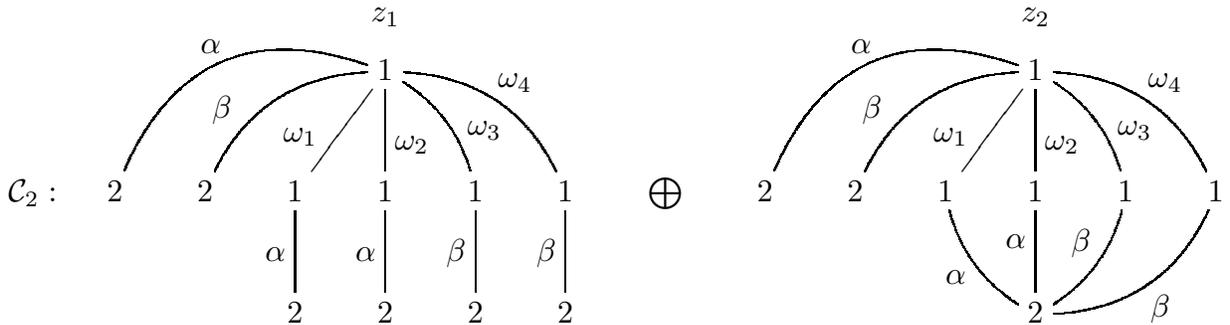
\xymatrix{
 &&&&1 \dropup{z_1} \edge@/_2pc/[dlll]_{\alpha} \edge@/_1pc/[dll]_(0.7){\beta} \edge[dl]_(0.65){\omega_1} \edge[d]^(0.6){\omega_2} \edge@/^/[dr]^(0.65){\omega_3} \edge@/^1pc/[drr]^{\omega_4}  && &&&&&1 \dropup{z_2} \edge@/_2pc/[dlll]_{\alpha}  \edge@/_1pc/[dll]_(0.7){\beta} \edge[dl]_(0.65){\omega_1} \edge[d]^(0.6){\omega_2} \edge@/^/[dr]^(0.65){\omega_3} \edge@/^1pc/[drr]^{\omega_4}  \\
\dropvert0{{\Cal C}_2:} &2 &2 &1 \edge[d]_{\alpha} &1 \edge[d]_{\alpha} &1 \edge[d]_{\beta} &1 \edge[d]_{\beta} &\bigoplus &2 &2 &1 \edge@/_/[dr]_{\alpha} &1 \edge[d]_(0.4){\alpha} &1 \edge@/^/[dl]_(0.45){\beta} &1 \edge@/^1pc/[dll]^{\beta}  \\
 &&&2 &2 &2 &2 & \dropvert{-7}{\bigl( c_1\alpha\omega_1 + c_2\beta\omega_4 \bigr) z_2 = 0  \text{\;\;for\;\;} [c_1:c_2] \in {\Cal C}_2} & &&&2
 }$$

\noindent {\it Partial reasoning\/}:  If we can show that the modules in $\C_1 \cup \C_2$ are degenerations of $M$, we know that they are maximal among those with top $T$;  indeed, from Theorem 6.5 it is immediate that the exhibited modules are devoid of proper top-stable degenerations.
We will only show how to recognize the modules in $\C_2$ as degenerations of $M$:  For $[c_1:c_2] \in \PP^1$,  consider the following curve $\phi: \PP^1 \rightarrow \overline{\autlap.C}$, determined by 
its restriction to $\AA^1$.  We define $\phi: \AA^1 \rightarrow \autlap.C$ so that, for $\tau \in \AA^1$, the map $\phi(\tau) \in \autlap$ sends $z_1$ to $z_1$ and $z_2$ to $z_2 + \tau c_1 \alpha \omega_1 z_2 +   \tau c_2 \beta \omega_4 z_2$.  It is readily checked that $\phi(\infty) = D$, which yields $M \degen P/D$ as claimed.    

That the union $\C_1 \cup \C_2$ includes all maximal top-stable degenerations of $M$ is immediate from Theorem 6.5:  Indeed, the annihilator in $M$ of the ideal $L$ has dimension $18$, and hence $\dim \ann_{M'} L \ge 18$ for any degeneration $M'$ of $M$.  Given that any maximal top-stable degeneration of $M$ is a direct sum of two local modules, at least one of the summands is therefore annihilated by $L$.  Now we invoke the additional restraints placed on maximal top-stable degenerations (Theorem 6.5(2)).   

Finally, we mention that the degenerations in $\C_1$ have height $3$ over $M$ (within the poset of degenerations of $M$), whereas those in $\C_2$ have height $2$.  On closer inspection, one moreover observes that the top-stable degenerations of $M$ of any fixed height above $M$ can be classified by a fine moduli space, each coming with an explicitly computable universal family.  The bulk of the top-stable degenerations of $M$ has height $1$.  It is a $\PP^3 \times \PP^1$-family which generically consists of indecomposable modules. 
\enddefinition

\head 7. Slicing $\lamod$ more finely, in terms of radical layerings   Representation-theoretically optimal coordinatization of $\grasstbd$ \endhead

In this section, we exhibit the features of $\grasstbd$ that provide representation-theoretic leverage.  Moreover, we will point to promising directions for uncovering further  classes of modules that permit classification through fine moduli spaces. The motto in doing so is to wield the knife in a manner guided by structural criteria.  Both objectives rest on an additional stratification of $\grasstbd$, namely into strata consisting of the points that represent modules with fixed radical layering.  (Caveat:  This partition of $\grasstbd$ into locally closed subvarieties is not necessarily a stratification in the strict technical sense, in that closures of strata need not be unions of strata in general.)  

Let $\SS = (\SS_0, \SS_1, \dots, \SS_L)$ be a sequence of semisimples in $\lamod$, where $L+1$ is the Loewy length of $\la$.  We call $\SS_0$ the {\it top\/} and $\underbardim \SS = \underbardim \bigoplus_{0 \le i \le L} \SS_l$ the {\it dimension vector\/} of $\SS$.  The sequences in which we are presently interested are the radical layerings $\SS(M) = (J^l M/J^{l+1} M)_{0 \le l \le L}$ of modules $M$.  Letting $\SS$ be a semisimple sequence with top $T$ and dimension vector $\bd$, we consider the following subvariety $\grassSS$ of $\grasstbd$:
$$\grassSS = \{C \in \grasstbd \mid \SS(P/C) = \SS\}.$$
Clearly, $\grassSS$ is stable under the $\autlap$-action on $\grasstbd$, where $P$ again denotes a fixed projective cover of $T$.  We will next introduce a representation-theoretically defined open affine cover of $\grassSS$, the charts of which are stable under the action of $\unirad$; in terms of stability of the charts, we can actually do a bit better.

In the first part of the upcoming definition, we describe a suitable basis for $P$.  For that purpose, we fix a sequence $z_1, \dots, z_t$ of top elements of $P$;  this means that the $z_r + JP$ form a basis for $P/JP$, and each $z_r$ is normed by one of the primitive idempotents, i.e., $z_r = e(r) z_r$ for some $e(r) \in \{e_1, \dots, e_n\}$.  Our choice clearly entails $P = \bigoplus_{1 \le r \le t} \la z_r$ with $\la z_r \cong \la e(r)$.

By a {\it path of length $l$ in\/} $P$ we mean a nonzero element of the form $p^{(l)} z_r$, where $p^{(l)}$ is a path of length $l$ in $KQ$.  Note that $p^{(l)}$ is then necessarily a path in $KQ \setminus I$ starting in $e(r)$.  The label $l$ serves to keep track of the length of the path in $KQ$, which is necessitated by the lack of an unambiguous concept of path length in $KQ/I$ in general.  Note that $z_1, \dots, z_t$ are precisely the paths of length $0$ in $P$.

\definition{Definition of skeleta and subsidiary comments}  
\smallskip

\noindent {\bf (I)}  A {\it skeleton of $P$\/} is a basis $B$ for $P$ with the following properties:

$\bullet$ $B$ consists of paths in $P$;

$\bullet$  For each $l \in \{0, \dots, L\}$, the cosets $p^{(l)} z_r + J^{l+1} P$ of the paths of length $l$ in $B$ form a basis for $J^l P/ J^{l+1} P$; 

$\bullet$ $B$ is closed under initial subpaths, that is: Whenever $p^{(l)} z_r \in B$ and $p^{(l)} = u_2^{(l_2)} u_1^{(l_1)}$ with paths $u_j^{(l_j)} \in KQ$, the path $u_1^{(l_1)} z_r  \in P$ belongs to $B$. 
\smallskip  

\noindent {\it Comments\/}:  $P$ has at least one skeleton, and any skeleton of $P$ contains $z_1, \dots, z_r$.  There is precisely one skeleton of $P$ in case $\la$ is a monomial algebra. Indeed, in the monomial case, the set of all $p_r^{(l)} z_r$, $1 \le r \le t$, where $p_r^{(l)}$ traces the paths of length $l$ in $KQ e(r) \setminus I$, is the unique skeleton, and the labels recording path lengths become superfluous.
\medskip

In the following, we fix a skeleton $B$ of $P$.
\smallskip

\noindent {\bf (II)}  An {\it {\rm{(}}abstract\/{\rm{)}} skeleton with radical layering $\SS = (\SS_0, \dots, \SS_L)$\/} is  
any subset $\S$ of $B$ which is closed under initial subpaths and is compatible with $\SS$ in the following sense:  For each $0 \le l \le L$, the multiplicity of $S_i$ in $\SS_l$ equals the number of those paths in $\S$ which have length $l$ and end in the vertex $e_i$.  
\smallskip

\noindent{\it Comment\/}: If $\SS_0 = T$, every skeleton with radical layering $\SS$ contains $z_1,\dots,z_t$. 
\smallskip

\noindent {\bf (III)}  Let $\S$ be an abstract skeleton with radical layering $\SS$.  We set
$$\grassS = \{C \in \grassSS \mid P/C \text{\ has basis\  } \S\}.$$
Whenever $M \in \lamod$ is isomorphic to some $P/C$ with $C \in \grassS$, we say that $\S$ is a {\it skeleton of $M$\/}.
\smallskip

\noindent{\it Comment\/}: Clearly, each $\grassSS$ is covered by the $\grassS$ that correspond to skeleta $\sigma$ with radical layering $\SS$ (finite in number). In other words, every module has at least one skeleton.
\smallskip
\enddefinition 

\definition{First consequences} Let $\S$ be any skeleton with radical layering $\SS$. 
\smallskip

 {\bf (1)} The set $\grassS$ is an open subvariety of $\grassSS$.  This is due to the following fact:  If $\SS$ has top $T$ and dimension $d$, then $\grassS$ is the intersection of $\grassSS$ with the big Schubert cell 
$$\Schu(\S) = \{\, C \in \Gr(\dim P - d,\, P) \mid P = C \oplus  \bigoplus_{b \in \S} Kb\, \}.$$  
Note, however, that $\grassS$ is not open in the ambient $\grasstbd$ (or $\grasstd$) in general. In particular, if $\SS$ has top $T$ and dimension vector $\bd$, the affine variety $\Schu(\S) \cap \grasstbd$ is typically larger than $\grassS$.

We infer that every irreducible component of $\grassSS$ comes with a generic set of skeleta.
\smallskip

{\bf (2)} Suppose $\S$ is a skeleton of $M$, i.e., $M = f(P/C)$ for some $C \in \grassS$ and isomorphism $f$. Then the radical layering of $M$ coincides with that of $\S$. Indeed, the definition entails that the paths of length $l$ in $\S$ induce a basis for the $l$-th radical layer of $M$, i.e., the residue classes $f(p^{(l)} z_r) + J^{l+1}M$, where the $p^{(l)} z_r$ run through the paths of length $l$ in $\S$, form a basis for $J^l M / J^{l+1} M$. 

For any choice $C \in \grassS$, we deduce: Whenever $b \in P \setminus \S$ is a path of length $l \ge 0$ and $\alpha$ is an arrow such that $\alpha b$ is a path in $P \setminus \S$, there exist unique scalars $c_{b'} \in K$ such that
$$\alpha b \equiv \sum_{b' \in \S(\alpha, b)} c_{b'} b' \pmod{C},$$
where $\S(\alpha, b)$ is the set of all paths in $\S$ which terminate in $\term(\alpha)$ and are at least as long as $b$. This places $\grassS$ into an affine space $\AA^N$ where $N$ is the sum of the cardinalities of the $\S(\alpha, b)$.
\smallskip

{\bf (3)} The variety $\grassS$ is affine. (Since the $\grassSS$ are not closed in $\grasstbd$, this requires proof; see \cite{\hier, Theorem 3.12}.) In fact, the coordinatization introduced in {\bf(2)} makes $\grassS$ a closed subset of $\AA^N$ up to isomorphism.  
\enddefinition

To make skeleta more user-friendly, we point to the fact that they are easy to visualize: any skeleton $\S$ may be identified with a forest.  Each tree in this forest consists of edge paths $p^{(l)} z_r$ for fixed $r$ and $l \ge 0$.  We illustrate this connection in a concrete situation.

\definition{Example 7.1}
Let
$\la = KQ/I$, where $Q$ is the quiver
$$\xymatrixcolsep{4pc}
 
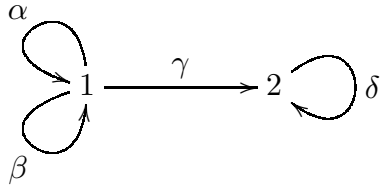
\xymatrix{ 
 1 \ar@'{@+{[0,0]+(0,8)}@+{[0,0]+(-8,8)}@+{[0,0]+(-8,2)}}_{\alpha} \ar@'{@+{[0,0]+(-8,-2)}@+{[0,0]+(-8,-8)}@+{[0,0]+(0,-8)}}_{\beta} \ar[r]^{\gamma} &2 \ar@'{@+{[0,0]+(6,6)}@+{[0,0]+(12,0)}@+{[0,0]+(6,-6)}}^{\delta}
 }$$
\noindent and $I$ is generated by all paths of length $4$.  Choose $T = S_1^3$, whence the projective cover of $T$ is $P = \bigoplus_{1 \le r \le 3} \la z_r$ with $\la z_r \cong \la e_1$ according to our convention.  Since $\la$ is a monomial algebra, $P$ has precisely one skeleton, namely the set 
$$B = \{p^{(l)} z_r \mid 1 \le r \le 3 \text{\;and\;} p^{(l)} \text{\ is a path of length\ } 0\le l \le 4 \text{\ in\ } KQ\setminus I\}.$$

Take $\bd = (6,5)$ and let $M \cong P/C$, where $C \in \grasstbd$ is generated by
$\gamma z_1$, $\gamma \alpha z_1$, $\alpha^2 z_1$, $\beta^2 z_1$,
 $(\beta \alpha - \alpha \beta)z_1$, $\alpha z_2$, $\beta z_2$, $\alpha z_3$, $\beta z_3$,
 $\gamma z_3 - \delta \gamma \beta z_1 - \delta \gamma z_2$. The point $C$ belongs to $\grassS$ for precisely two skeleta $\S$, which we present graphically below. (The alternate skeleta of $M$ result from permutations of the trees in the pertinent $\S$.)
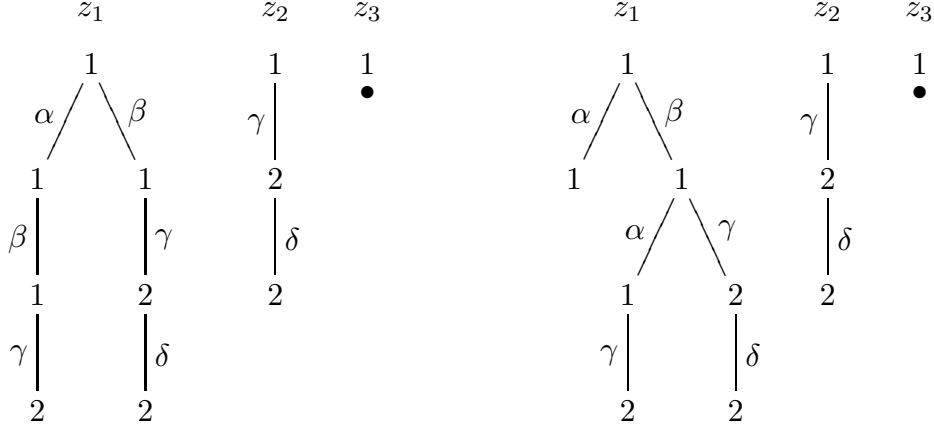
$$\xymatrixrowsep{2pc}\xymatrixcolsep{0.5pc}
\xymatrix{
 &1 \dropup{z_1} \edge[dl]_{\alpha} \edge[dr]^{\beta} &&&&1 \dropup{z_2}
\edge[d]_{\gamma} &&1 {\save+<0ex,-2ex> \drop{\bullet} \restore} \dropup{z_3} 
&&&&& &1 \dropup{z_1} \edge[dl]_{\alpha} \edge[dr]^{\beta} &&&&1 \dropup{z_2}
\edge[d]_{\gamma} &&1 {\save+<0ex,-2ex> \drop{\bullet} \restore} \dropup{z_3} \\
1 \edge[d]_{\beta} &&1 \edge[d]^{\gamma} &&&2 \edge[d]^{\delta} 
&&&&&&&1 &&1 \edge[dl]_{\alpha} \edge[dr]^{\gamma} &&&2 \edge[d]^{\delta} \\
1 \edge[d]_{\gamma} &&2 \edge[d]^{\delta}  &&&2
&&&&&&& &1 \edge[d]_{\gamma} &&2 \edge[d]^{\delta} &&2  \\
 2 &&2  &&&&&&&&&& &2 &&2 
}$$
\noindent In formal terms, the first skeleton consists of all the paths $p^{(l)} z_r$ that
occur as edge paths (of length $l \ge 0$) in one of the three left-hand trees as one reads them from top to
bottom. 

The radical layering of $M$ is equal to the radical layerings of the above skeleta, namely $\SS(M) = (S_1^3, \, S_1^2 \oplus S_2, \, S_1 \oplus S_2^2, \, S_2^2)$.
\enddefinition

Consequence {\bf(2)} above guarantees that the final observation we made in the example generalizes: From any skeleton $\S$ of a module $M$, we retrieve the radical layering of $M$. 
\medskip

We return to the general discussion, to address stability properties of the $\grassS$.
Based on our choice of top elements $z_1, \dots, z_t$ of $P$, we pin down a maximal torus in $\autlap$:  Namely, we let  $\T$ be  the group of automorphisms $P \rightarrow P$ defined by $z_r \mapsto a_r z_r$, for some element $(a_1, \dots, a_t)$ in the torus $(K^*)^t$.

One of the crucial levers, applied (e.g.) to the proof of the classification results in Section 6, is as follows.  It rests on theorems due to Kostant and Rosenlicht (see \cite{\Rosone, Theorem 2} and \cite{\Rostwo, Theorem 1}):   Any morphic action of a unipotent group on an affine variety has closed orbits, and these orbits are full affine spaces.  The former fact clearly entails the second of the following assets of the affine cover $\bigl( \grassS \bigr) _\S$. 

\proclaim{Proposition 7.2}  Let $T$ be an arbitrary semisimple module, and $P$ be as before.
\smallskip 

$\bullet$  For every skeleton $\S$ with top $T$, the variety $\grassS$ is stable under the action of 

$\T \ltimes \unirad$.

In particular:  If $T$ is squarefree, the $\grassS$ are stable under $\autlap$.
\smallskip

$\bullet$ For every semisimple sequence $\SS$ with top $T$, the $\unirad$-orbits are closed in 

$\grassSS$.  

In particular:  If $T$ is squarefree, all $\autlap$-orbits of $\grassSS$ are closed in 
$\grassSS$.
\endproclaim

The final statement of the proposition (as well as the simple structure of the $\unirad$-orbits) explains why modules with squarefree tops hold a special place in the exploration of orbit closures.
Unfortunately, relative closedness of the $\autlap$-orbits of $\grassSS$ does not, by itself, guarantee existence of an orbit space of $\grassSS$ by $\autlap$.  Another necessary condition is that $\autlap$ act with constant orbit dimension on the irreducible components of $\grassSS$ (see, e.g., \cite{\Bor, Chapter II, Proposition 6.4}).  This orbit-equidimensionality commonly fails;  discrepancies among the orbit dimensions on irreducible components of $\grassSS$ may actually be arbitrarily large.  However, the situation can often be salvaged through some additional slicing guided by skeleta.   We include two examples to illustrate this tack at classification.  These instances are by no means isolated  --   in fact, we cannot name an example of a sequence $\SS$ with squarefree top where the underlying technique fails  --   but a systematic investigation along this line has not been undertaken so far. 

In 7.3 and 7.4, we specify choices of $\la$, $T$, and $\bd$. Each time, we will encounter the following situation: The modules with top $T$ and dimension vector $\bd$ do not have a coarse moduli space.  Yet, for any radical layering $\SS$ with that top and dimension vector, the modules in $\grassSS$ either have a fine moduli space, or else have a finite partition into locally closed subsets, specifiable in terms of module structure, such that the representations parametrized by the individual subvarieties are finely classifiable.

\definition{Example 7.3}  Let $\la = KQ/I$, where  $Q$ is the quiver 
$$\xymatrixrowsep{3pc}\xymatrixcolsep{6pc}
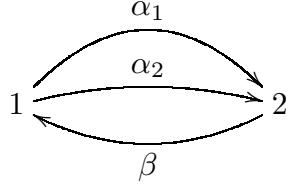
\xymatrix{
1 \ar@/^2pc/[r]^{\alpha_1} \ar@/^/[r]^{\alpha_2}  &2  \ar@/^1pc/[l]^{\beta}
}$$
\noindent and $I$ is the ideal  generated by all paths of length $3$. Moreover, we choose $T = S_1 \oplus S_2$ and $\bd = (2,2)$.  In accordance with our conventions, we write $P = \la z_1 \oplus \la z_2$ with $z_j = e_j$. By Corollary 6.6, the modules with top $T$ and dimension vector $\bd$ do not have a fine moduli space; indeed, the point $C =  \la \alpha_2 z_1 +  \la \alpha_1 \beta z_2 + \la (\alpha_1 z_1 - \alpha_2 \beta z_2)$ in $\grasstbd$, for instance, is not invariant under automorphisms of $P$.  

We analyze the three strata $\grassSS$ that make up $\grasstbd$. They are all irreducible, by \cite{\BHTtwo, Theorem 5.3}.  For $\SS = (T, T, 0)$, the modules in $\grassSS$ clearly have a fine moduli space, namely $\PP^1$.  
For $\SS = (T, S_2, S_1)$, the outcome is the same.

Now we focus on $\SS = (T, S_1, S_2)$. The orbit dimension being non-constant on this stratum, we subdivide it further to arrive at classifiable portions. The class $\C_1$ of decomposable modules with radical layering $\SS$ is easily seen to have a fine moduli space, namely $\PP^1$; the class $\C_2$ of indecomposable modules with skeleton $\S = \{z_1, z_2, \alpha_1 z_1, \beta z_2, \alpha_2 \beta z_2\}$ has a fine moduli space as well, namely $\AA^1$.  To verify this, observe that each module in $\C_2$ has a normal form $P/C_k$, where 
$$C_k  = \la \alpha_2 z_1 + \la (  \alpha_1 z_1 - \alpha_2 \beta z_2) + \la ( \alpha_1 \beta - k \alpha_2 \beta) z_2$$ 
for a unique scalar $k \in K$; the universal family for $\C_2$ is the corresponding trivial bundle.  There is only a single $\autlap$-orbit in $\grassSS$ which does not belong to $\C_1 \cup \C_2$, namely that of $C = \la \alpha_1 z_1 +  \la \alpha_2 \beta z_2 + \la \bigl(  \alpha_2 z_1 - \alpha_1 \beta z_2 \bigr)$. 
\enddefinition

\definition{Example 7.4}  This time, start with the algebra $\la = KQ/I$, where $Q$ is the quiver 
$$\xymatrixcolsep{4.5pc}
 
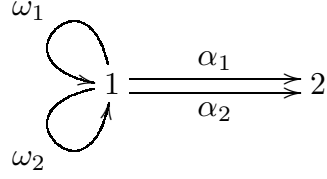
\xymatrix{ 
 1 \ar@'{@+{[0,0]+(-1,8)}@+{[0,0]+(-8,8)}@+{[0,0]+(-8,1)}}_{\omega_1} \ar@'{@+{[0,0]+(-8,-1)}@+{[0,0]+(-8,-8)}@+{[0,0]+(-1,-8)}}_{\omega_2} 
  \ar[r]<0.5ex>^{\alpha_1}  \ar[r]<-0.5ex>_{\alpha_2} &2
 }$$
\noindent  and $I$ is generated by all $\omega_i \omega_j$ together with the paths $\alpha_i \omega_j$ for $i \ne j$.  Let $T = S_1$ (hence $P = \la z_1$ with $z_1 = e_1$), and $\bd = (3, 1)$.

If $\SS = (S_1, S_1^2 \oplus S_2, 0)$, it is easy to see that the isomorphism classes of modules with radical layering $\SS$ have a fine moduli space, namely $\PP^1$.

 Now suppose that  $\SS = (S_1, S_1^2, S_2)$, the only alternate radical layering with top $T$ and dimension vector $\bd$. Once again, the variety $\grassSS$ is irreducible, but the dimensions of its $\autlap$-orbits fail to be constant. Hence, the modules with radical layering $\SS$ do not have a fine moduli space.

On the other hand, let $\C_1$ be the class of modules $M$ with radical layering $\SS$ such that $\alpha_j \omega_j M \ne 0$ for $j=1,2$, and let $\C_2$ consist of the remaining modules with radical layering $\SS$.  Each of $\C_1$, $\C_2$ has a fine moduli space providing classification up to isomorphism.

Indeed, the modules in $\C_1$ have graphs of the form 
$$\xymatrixrowsep{2pc}\xymatrixcolsep{0.5pc}
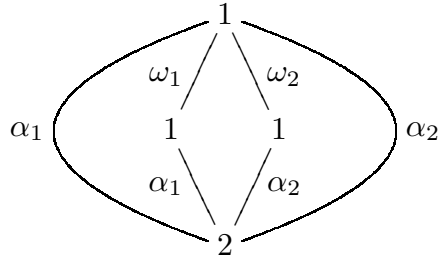
\xymatrix{
 &1 \edge[dl]_(0.6){\omega_1} \edge@/_4.5pc/[dd]_{\alpha_1} \edge[dr]^(0.6){\omega_2} \edge@/^4.5pc/[dd]^{\alpha_2}  \\
1 \edge[dr]_(0.4){\alpha_1} &&1 \edge[dl]^(0.4){\alpha_2}  \\
 &2
}$$
\noindent and the set $X_1$ of all points in $\grassSS$ corresponding to modules in $\C_1$ is an $\autlap$-stable subvariety isomorphic to $K^* \times \AA^2$.  It is readily checked that, for any module $M$ in $\C_1$, there exists a unique scalar $k \in K^*$ such that $M \cong P / C_k$ with $C_k = \la (\alpha_2 \omega_2 - k \alpha_1 \omega_1) z_1 + \sum_{j = 1,2} \la \alpha_j z_1$. In fact, the canonical projection $X_1 = K^* \times \AA^2  \rightarrow K^*$ is a geometric quotient of $X_1$ by $\autlap$.  The trivial bundle $\Delta = K^* \times K^4$, endowed with the $K$-algebra homomorphism $\delta: \la \rightarrow \End(\Delta)$ that is obtained along the preceding recipe, is the corresponding universal family, confirming that $K^*$ is a fine moduli space for $\C_1$.

Analogously, the modules in $\C_2$, represented by the orbits in $X_2 = \grassSS \setminus X_1$, have graphs
$$\xymatrixrowsep{2pc}\xymatrixcolsep{1.75pc}
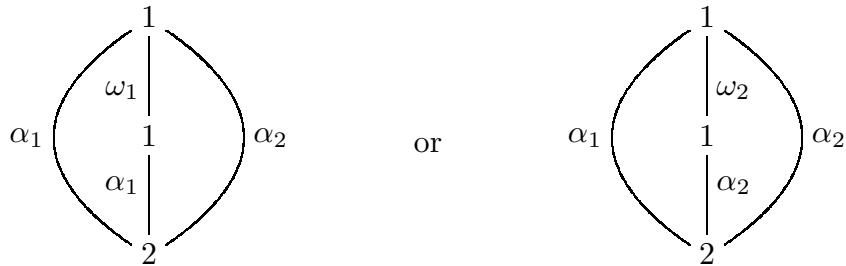
\xymatrix{
1 \edge[d]_(0.6){\omega_1}  \edge@/_2.5pc/[dd]_{\alpha_1} \edge@/^2.5pc/[dd]^{\alpha_2}  &&&&&&1  \edge[d]^(0.6){\omega_2}  \edge@/_2.5pc/[dd]_{\alpha_1} \edge@/^2.5pc/[dd]^{\alpha_2}  \\
1 \edge[d]_(0.4){\alpha_1} &&&\txt{or} &&&1 \edge[d]^(0.4){\alpha_2}  \\  
2 &&&&&&2
 }$$
\noindent  depending on whether $\alpha_1 \omega_1$ or $\alpha_2 \omega_2$ annihilates.  In fact, $X_2$ consists of two disjoint irreducible components, reflecting the dichotomy with respect to annihilators.  Considerations following the previous pattern yield normal forms 
$P/C_k$, with $C_k = \la( \alpha_2 - k \alpha_1 \omega_1) z_1 + \la \alpha_1 z_1 + \la \alpha_2 \omega_2 z_1$ for $k \in K$, for the modules in $\C_2$ with skeleton $\{ z_1, \omega_1 z_1, \alpha_1 \omega_1 z_1 \}$. A symmetric description applies to the modules in the other component of $X_2$.  Here, $k = 0$ is allowed as well.  For $k \ne 0$, the normal forms are reflected by simplified graphs 
$$\xymatrixrowsep{2pc}\xymatrixcolsep{1.75pc}
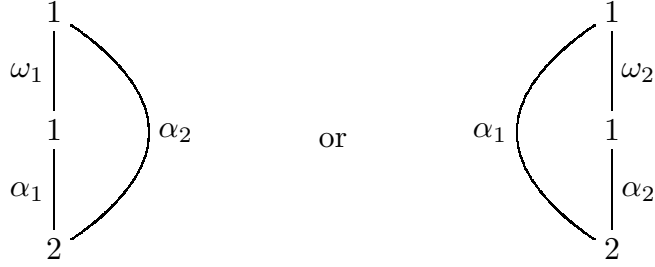
\xymatrix{
1 \edge[d]_{\omega_1}   \edge@/^2.5pc/[dd]^{\alpha_2}  &&&&&&1 \edge[d]^{\omega_2}  \edge@/_2.5pc/[dd]_{\alpha_1}  \\
1 \edge[d]_{\alpha_1}  &&&\txt{or} &&&1 \edge[d]^{\alpha_2}  \\
2 &&&&&&2
}$$
\noindent  Guided by this observation, one verifies that the class $\C_2$ has a fine moduli space as well, this one consisting of two irreducible components isomorphic to $\AA^1$. 
\enddefinition

In extending the ideas illustrated in 7.3 and 7.4 to semisimple modules $T = \bigoplus_{1 \le i \le n} S_i^{t_i}$, where multiplicities $\ge 2$ are permitted, the following straightforward observation turns out useful.

\proclaim{Observation 7.5}  Let $\sigma$ be a skeleton with radical layering $\SS$ and, again, let $P$ be a fixed projective cover of $T$.  The orbit-closure of any open subvariety $X$ of $\grassS$ is an open subvariety of $\grassSS$. 
 \endproclaim
 
 Such orbit-closures make the methods of the examples applicable to situations where the $\grassS$ fail to be $\autlap$-stable.

We conclude this section by placing a spotlight on the need to focus on non-closed subvarieties of the $\grasstbd$ to make broader use of the geometric classification tools we have described.  Suppose $T_1, \dots, T_m$ is a sequence of semisimple $\la$-modules, and consider the requirement that the modules with top $T_j$ be classifiable via moduli spaces for all $j$.  The (proof of) the upcoming equivalence demonstrates how the pressure placed on the algebra $\la$ rapidly builds as we enlarge the collection of $T_j$.
 
\proclaim{Proposition 7.6}  The following conditions on $\la$ are equivalent:
\smallskip

{\rm {\bf (a)}}  For all semisimple modules $T \in \lamod$ of dimension $2$ and any choice of $\bd$, the left $\la$-modules with top $T$ and dimension vector $\bd$ have a fine {\rm{(}}equivalently, a coarse{\rm{)}} moduli space classifying them up to isomorphism.
\smallskip

{\rm {\bf (b)}}  $\la$ is a Nakayama algebra, that is, all $\la$-modules are direct sums of uniserials.
\endproclaim

\demo{Proof}  In light of Theorems 6.3 and 6.5 it is immediate that (b) implies (a). (In fact, this implication does not require restriction to dimension $2$.)

For the reverse implication, assume (a).  Since (b) is equivalent to the requirement that all indecomposable projective left or right $\la$-modules be uniserial, we only need to show that any vertex in $Q$ is subject to the following constraint:  it does not occur as the starting point of more than one arrow, nor as the end point of more than one arrow.

First suppose that there is a vertex $e$ such that two distinct arrows $\alpha$ and $\beta$ start in $e$.  Let $T = (\la e/ Je)^2$ and $\bd$ the dimension vector of $T \oplus S \oplus S'$, where the final two summands are the simple modules corresponding to the terminal vertices of $\alpha$ and $\beta$, respectively (in particular, the possibility $S \cong S'$ is not excluded).      
Moreover, consider the following point $C = C_1 \oplus C_2 \in \grasstbd$:  Let $P = \la z_1 \oplus \la z_2$ where $\la z_1 \cong \la z_2 \cong \la e$;  define $C_1$ to be the submodule of  $\la z_1$ generated by all elements $p z_1$ where $p$ traces the arrows different from $\alpha$ and all paths of length $2$; define the submodule $C_2 \subseteq \la z_2$ symmetrically with $\beta$ taking over the role of $\alpha$.  Then $P/C$ has a proper top-stable degeneration by Theorem 6.5, since the left ideals of $\la e$ corresponding to $C_1$ and $C_2$ are not comparable.  Hence, the modules with top $T$ and dimension vector $\bd$ fail to have a coarse moduli space under the present assumption.

Now suppose that there is a vertex $e$ with two distinct arrows $\alpha$ and $\beta$ ending in $e$. Let $S$ and $S'$ be the simple modules corresponding to the starting vertices of $\alpha$ and $\beta$, respectively. Set $T = S \oplus S'$ and $\bd = \underbardim S \oplus S' \oplus (\la e / Je)$.  Again fix a projective cover of $T$, say $P = \la z \oplus \la z'$, where $z$ and $z'$ correspond to the starting points of $\alpha$ and $\beta$, respectively.  This time, $C \in \grasstbd$ is to be generated by $\alpha z - \beta z'$, all elements $p z$, $p z'$ where $p$ traces the paths of length $2$, next to all elements $p z$ where $p$ is an arrow different from $\alpha$, and all $p z'$ where $p$ is an arrow different from $\beta$.  The module $P/C$ is clearly indecomposable non-local, and hence again has a proper top-stable degeneration by Theorem 6.5.  Once more, this precludes existence of a coarse moduli space classifying the modules with top $T$ and dimension vector $\bd$. $\qed$
\enddemo

\head 8.  Problems.  Pros and Cons of Approach B \endhead

\definition{Open Problems}
\smallskip

The first series of problems consists of immediate followups to the results of Section 6.
\smallskip

(1) We saw that arbitrary projective varieties arise as fine moduli spaces $\maxmoduli$ for suitable choices of $\la$, $T$, and $\bd$.  This begs the question:  Which projective varieties arise as fine moduli spaces, $\maxmoduliM(M)$, classifying the maximal top-stable degenerations of an individual module $M$?
\smallskip

(2)  Relate the structure of $\maxmoduli$ (resp\. of $\maxmoduliM(M)$) to $\la$, $T$, and $\bd$ (resp\. to $M$).  In particular, investigate rationality (which typically facilitates the analysis of the generic structure of the modules in the irreducible components; see \cite{\BHTtwo}) and normality.
\smallskip

(3)  Let $M \in \lamod$.  In all presently known examples,  the top-stable degenerations of $M$ of fixed height above $M$ (in the poset of degenerations of $M$) have representation-theoretically defined finite partitions with the property that each of the corresponding isomorphism classes of degenerations of $M$ has a fine moduli space.  Explore this phenomenon systematically, beginning with the case of a simple top. 
\smallskip

(4)  Let $T'$ be a semisimple module properly containing the top $T$ of $M$.  Compare $\maxmoduliM(M)$ with the fine moduli space classifying the degenerations of $M$ which are maximal among those with top $T'$. 
\smallskip

The second set of problems is motivated by the observations and examples in Section 7.    
\smallskip

(5) Let $\SS$ be a sequence of semisimple modules with squarefree top.  
Confirm or refute the following equivalences:
 \smallskip

 $\bullet$  There is a coarse moduli space classifying the modules with radical layering $\SS$ up to isomorphism.

 $\bullet$  The $\autlap$-orbits (= $\unirad$-orbits) of $\grassSS$ have constant dimension.
\smallskip

(6)  For a given algebra $\la$ (from a specified class), determine the radical layerings $\SS$ such that the modules degeneration-maximal among those with radical layering $\SS$ have a coarse or fine moduli space. 

(Note:  In Example 7.4, all modules with radical layering $\SS = (S_1, S_1^2, S_2)$ are degenera\-tion-maximal among those with the given radical layering, since, in the case of a simple top, $\autlap$ operates with closed orbits on $\grassSS$; see Proposition 7.2.  However, $\grassSS$ does not possess an orbit space.) 

\smallskip

(7)  Develop a general slicing technique for the varieties $\grassSS$ on the model of Examples 7.3 and 7.4, at least for sequences $\SS$ with squarefree tops.  

(The idea is to find partitions that are not ``opaque" from a representation-theoretic viewpoint, so  as not to defeat the purpose of classification.)
\smallskip

(8)  Generically classify the modules in $\grassSS$ in the following sense:  For each irreducible component Comp of $\grassSS$, specify, in representation-theoretic terms, a dense open subvariety $X$(Comp) such that the modules parametrized by $X$(Comp) have a fine or coarse moduli space. 

(Note:  In light of \cite{\BHTtwo}, the irreducible components of arbitrary varieties $\grassSS$ can be algorithmically determined from a presentation of the underlying algebra in terms of quiver and relations.  In case $\la$ is a truncated path algebra, all of the pertinent $\grassSS$ are irreducible, whence this special case provides a good starting point.

Nontrivial instances of generic classification: In Example 7.3, the variety $\grassSS$ for $\SS = (S_1 \oplus S_2, S_1, S_2)$ is irreducible, and the subvariety consisting of the orbits that correspond to the modules in $C_2$ yields a generic classification.  In Example 7.4, take $\SS = (S_1, S_1^2, S_2)$, which again leads to an irreducible variety $\grassSS$; this time the $\autlap$-stable subvariety corresponding to the modules in the class $\C_1$ provides a generic classification of the desired ilk.)

\enddefinition
\medskip

\centerline {{\bf Pros and cons of Approach B}}
\smallskip
\noindent  {\bf Pros:}

$\bullet$  One controls the class of modules to be classified.  In fact, the target classes are cut out of $\lamod$ in representation-theoretic terms to begin with.  The same is true for the equivalence relation up to which one is trying to classify:  It is either isomorphism, or else isomorphism preserving some additional structure (such as a grading).   
\smallskip

$\bullet$  In the instances addressed in Section 6, the moduli spaces are quite accessible to computation, in a manner that ties their geometry to the combinatorics of a presentation by quiver and relations of the underlying algebra.  The same holds for the construction of the corresponding universal families.  (The computational access is via the closed embeddings of these moduli spaces into the pertinent $\grasstbd$ and the computable affine charts $\grassS$;  an algorithm for finding the $\grassS$ from the quiver $Q$ and generators of $I$ has been implemented.  Due to the transparent connection between points of $\grasstbd$ and minimal projective presentations of the modules they encode, constructing restrictions of the universal families to the affine charts follows suit.) 
\medskip

\noindent {\bf Cons:}

$\bullet$  Existence of coarse or fine moduli spaces for the representations corresponding to large closed subvarieties of the $\grasstbd$ is a rare occurrence (see Proposition 7.6 and the comments following Theorem 6.5).   There is no machinery that guarantees existence of moduli spaces coming out of specific search strategies.  Here, in turn, there is considerable reliance on serendipity, just different in nature from that required in Approach A.  Under the latter strategy, one relies on effective choices of weight functions, while under Strategy B one relies on the discovery of promising normal forms of the target classes of representations.  For neither task is there a general recipe.
\smallskip 

$\bullet$  As one moves beyond the instances of classifiability exhibited in Section 6, one is likely to sacrifice grasp of most geometric boundary phenomena arising in the varieties $\grasstbd$ (once again, see Section 7).  This downside parallels one of the negatives singled out in connection with Approach A.
\smallskip    

$\bullet$  There is no ``ready-made" arsenal of techniques available for the geometric analysis of the resulting moduli spaces, existence provided.  Followup methods for taking optimal advantage of existence results need to be designed to measure, case by case.

\enddocument